\documentclass{aims}
\usepackage{amsmath}
  \usepackage{paralist}
  \usepackage{graphics} 
  \usepackage{epsfig} 
\usepackage{graphicx} 
\usepackage{slashbox,booktabs,amsmath} 
\usepackage{epstopdf}
\usepackage[colorlinks=true]{hyperref}
\hypersetup{urlcolor=blue, citecolor=red}

\def\currentvolume{}
 \def\currentissue{}
  \def\currentyear{}
   \def\currentmonth{}
    \def\ppages{}
     \def\DOI{}

\newtheorem{theorem}{Theorem}[section]

\newtheorem{lemma}[theorem]{Lemma}

\theoremstyle{definition}
\newtheorem{definition}[theorem]{Definition}

\title[Polynomial optimization for stability analysis and Control] 
      {Polynomial optimization with applications to stability analysis and Control - Alternatives to Sum of Squares}

\author[Reza Kamyar and Matthew M. Peet]{}


 \email{rkamyar@asu.edu}
 \email{mpeet@asu.edu}

\thanks{Acknowledgements: NSF grant No. CMMI-1151018 and CMMI-1100376}

\begin{document}
\maketitle

\centerline{\scshape Reza Kamyar }
\medskip
{\footnotesize
   \centerline{School of Matter, Transport and Energy}
   \centerline{Arizona State University}
   \centerline{Goldwater Center for Science and Engineering 531}
   \centerline{Tempe, 85281, USA}
} 

\medskip

\centerline{\scshape Matthew M. Peet}
\medskip
{\footnotesize
   \centerline{School of Matter, Transport and Energy}
   \centerline{Arizona State University}
   \centerline{501 Tyler Mall - ECG 301}
   \centerline{Tempe, 85281, USA}
}

\bigskip

 \centerline{}

\begin{abstract}
In this paper, we explore the merits of various algorithms for polynomial optimization problems, focusing on alternatives to sum of squares programming. While we refer to  advantages and disadvantages of Quantifier Elimination, Reformulation Linear Techniques, Blossoming and Groebner basis methods, our main focus is on algorithms defined by Polya's theorem, Bernstein's theorem and Handelman's theorem. We first formulate polynomial optimization problems as verifying the feasibility of semi-algebraic sets. Then, we discuss how Polya's algorithm, Bernstein's algorithm and Handelman's algorithm reduce the intractable problem of feasibility of semi-algebraic sets to linear and/or semi-definite programming. We apply these algorithms to different problems in robust stability analysis and stability of nonlinear dynamical systems. As one contribution of this paper, we apply Polya's algorithm to the problem of $H_\infty$ control of systems with parametric uncertainty.
 Numerical examples are provided to compare the accuracy of these algorithms with other polynomial optimization algorithms in the literature.

\end{abstract}

\section{Introduction}
Consider problems such as portfolio optimization, structural design, local stability of nonlinear ordinary differential equations, control of time-delay systems and control of systems with uncertainties. These problems can all be formulated as \textit{polynomial optimization} or \textit{optimization of polynomials}. In this paper, we survey how computation can be applied to polynomial optimization and optimization of polynomials. 
One example of polynomial optimization is $\beta^* = \min_{x \in \mathbb{R}^n} p(x)$, where $p: \mathbb{R}^n \rightarrow \mathbb{R}$ is a multi-variate polynomial. In general, since $p(x)$ is not convex, this is not a convex optimization problem. It is well-known that polynomial optimization is NP-hard~\cite{blum}. 
We refer to the dual problem to polynomial optimization as optimization of polynomials, e.g., the dual optimization of polynomials to $\beta^* = \min_{x \in \mathbb{R}^n} p(x)$ is 
\begin{align}
& \quad \beta^* =\max_{y \in \mathbb{R}} \; y \nonumber \\
&\text{subject to } \;\; p(x) - y \geq 0 \text{ for all } x \in \mathbb{R}^n.
\label{optim_poly_exp}
\end{align}
This problem is convex, yet NP-hard. 

One approach to find lower bounds on the optimal objective $\beta^*$ is to apply Sum of Squares (SOS) \textit{programming}~\cite{parillo_thesis,sostools2013}. A polynomial $p$ is SOS if there exist polynomials $q_i$ such that $p(x)=\sum_{i=1}^r q_i(x)^2$.
The set $\{ q_i \in \mathbb{R}[x], i=1, \cdots,r \} $ is called an SOS \textit{decomposition} of $p(x)$, where $\mathbb{R}[x]$ is the ring of real polynomials. An SOS program is an optimization problem of the form
\begin{align}
& \quad \min_{x \in \mathbb{R}^m} \qquad c^T x \nonumber \\
&\text{subject to} \quad  A_{i,0}(y)+\sum_{j=1}^m x_j A_{i,j}(y) \text{ is SOS, } i=1, \cdots, k,
\label{eq:SOS_prog}
\end{align}
where $c \in \mathbb{R}^m$ and $A_{i,j} \in \mathbb{R}[y]$ are given. If $p(x)$ is SOS, then clearly $p(x) \geq 0$ on $\mathbb{R}^n$. While verifying $p(x) \geq 0$ on $\mathbb{R}^n$ is NP-hard, checking whether $p(x)$ is SOS - hence non-negative - can be done in polynomial time~\cite{parillo_thesis}. It was first shown in~\cite{parillo_thesis} that verifying the existence of a SOS decomposition is a Semi-Definite Program. Fortunately, there exist several algorithms~\cite{monteiro, helmberg, alizadeh} and solvers~\cite{sdpa,sedumi,sdpt3} that solve SDPs to arbitrary precision in polynomial time. To find lower bounds on $\beta^* = \min_{x \in \mathbb{R}^n} p(x)$, consider the SOS program
\[
y^*=\max_{y \in \mathbb{R}}  y \; \text{ subject to } \; p(x) - y \text{ is SOS}. 
\]
Clearly $y^* \leq \beta^*$. By performing a bisection search on $y$ and semi-definite programming to verify $p(x) - y$ is SOS, one can find $y^*$. SOS programming can also be used to find lower bounds on the global minimum of polynomials over a semi-algebraic set $S:=\{x \in \mathbb{R}^n: g_i(x) \geq 0, h_j(x)=0 \}$ generated by $g_i,h_j \in \mathbb{R}[x]$.
Given problem~\eqref{optim_poly_exp} with $x \in S$, \textit{Positivstellensatz} results~\cite{stengle,putinar,schmudgen} define a sequence of SOS programs whose objective values form a sequence of lower bounds on the global minimum $\beta^*$. It is shown that under certain conditions on $S$~\cite{putinar}, the sequence of lower bounds converges to the global minimum. See~\cite{laurent} for a comprehensive discussion on the Positivstellensatz. 

In this paper, we explore the merits of some of the alternatives to SOS programming. There exist several results in the literature that can be applied to polynomial optimization; e.g., Quantifier Elimination (QE) algorithms~\cite{CAD} for testing the feasibility of semi-algebraic sets, Reformulation Linear Techniques (RLTs)~\cite{sherali_1992,sherali_1997} for linearizing polynomial optimizations, Polya's result~\cite{inequalities} for positivity on the positive orthant, Bernstein's~\cite{roy,leroy} and Handelman's~\cite{handelman_1988} results for positivity on simplices and convex polytopes, and other results based on Groebner bases~\cite{adams_groebner} and Blossoming~\cite{blossoming}. We will discuss Polya's, Bernstein's and Handelman's results in more depth. The discussion of the other results are beyond the scope of this paper, however the ideas behind these results can be summarized as follows. 

QE algorithms apply to First-Order Logic formulae, e.g., 
\[
\forall x \, \exists y \, (f(x,y) \geq 0 \Rightarrow ((g(a) < x y) \wedge (a >  2)),
\]
to eliminate the \textit{quantified} variables $x$ and $y$ (preceded by quantifiers $\forall,\exists$) and construct an equivalent formula in terms of the \textit{unquantified} variable $a$. The key result underlying QE algorithms is Tarski-Seidenberg theorem~\cite{tarski}. The theorem implies that for every formula of the form $\forall x \in \mathbb{R}^n \, \exists y \in \mathbb{R}^m ( f_i(x,y,a) \geq 0 )$, where $f_i \in \mathbb{R}[x,y,a]$, there exists an equivalent quantifier-free formula of the form $\wedge_i (g_i(a) \geq 0) \vee_j (h_j(a) \geq 0)$ with $g_i,h_j \in \mathbb{R}[a]$. QE implementations~\cite{QEPCAD,Redlog} with a bisection search yields the exact solution to optimization of polynomials, however the complexity scales double exponentially in the dimension of variables $x,y$.

RLT was initially developed to find the convex hull of feasible solutions of zero-one linear programs~\cite{sherali_1990}. It was later generalized to address polynomial optimizations of the form $\min_x p(x)$ subject to $x \in [0,1]^n \cap S$~\cite{sherali_1992}. RLT constructs a $\delta-$hierarchy of linear programs by performing two steps. In the first step (reformulation), RLT introduces the new constraints $\prod_i x_i \prod_j (1-x_j) \geq 0$ for all $i,j: i+j=\delta$. In the second step (linearization), RTL defines a linear program by replacing every product of variables $x_i$ by a new variable. By increasing $\delta$ and repeating the two steps, one can construct a $\delta-$hierarchy of lower bounding linear programs. A combination of RLT and branch-and-bound partitioning of $[0,1]^n$ was developed in~\cite{sherali_1997} to achieve tighter lower bounds on the global minimum. For a survey of different extensions of RLT see~\cite{sherali_global_2007}. 

Groebner bases can be used to reduce a polynomial optimization over a semi-algebraic set $S:=\{ x \in \mathbb{R}^n : g_i(x) \geq 0, \, h_j(x) = 0\}$ to the problem of finding the roots of univariate polynomials~\cite{PP_groebner}. First, one needs to construct the system of polynomial equations
\begin{equation}
\nabla_x L(x,\lambda,\mu), \nabla_\lambda L(x,\lambda,\mu), \nabla_\mu L(x,\lambda,\mu)] = 0,
\label{eq:nabla_sys}
\end{equation}
 where $L:= p(x)+ \sum_i \lambda_i g_i(x)+ \sum_j \mu_j h_j(x)$ is the Lagrangian.
 It is well-known that the set of solutions to~\eqref{eq:nabla_sys} is the set of extrema of the polynomial optimization $\min_{x \in S} p(x)$.
Let
 \[
\left[ f_1(x,\lambda,\mu), \cdots,f_N(x,\lambda,\mu) \right] := \left[\nabla_x L(x,\lambda,\mu), \nabla_\lambda L(x,\lambda,\mu), \nabla_\mu L(x,\lambda,\mu) \right].
 \]
Using the elimination property~\cite{adams_groebner} of the Groebner bases, the minimal Groebner basis of the ideal of $f_1, \cdots, f_N$ defines a triangular-form system of polynomial equations. This system can be solved by calculating one variable at a time and back-substituting into other polynomials. The most computationally expensive part is the calculation of the Groebner basis, which in the worst case scales double-exponentially in the number of decision variables.

The blossoming approach involves mapping the space of polynomials to the space of multi-affine functions (polynomials that are affine in each variable). By using this map and the diagonal property of blossoms~\cite{blossoming}, one can reformulate any polynomial optimization $\min_{x \in S} p(x)$ as an optimization of multi-affine functions. In~\cite{girard_blossom}, it is shown that the dual to optimization of multi-affine functions over a hypercube is a linear program. The optimal objective value of this linear program is a lower bound on the minimum of $p(x)$ over the hypercube.

While the discussed algorithms have advantages and disadvantages (such as exponential complexity), we focus on Polya's, Bernstein's and Handelman's theorems - results which yield polynomial-time parameterizations of positive polynomials. Polya's theorem yields a basis to parameterize the cone of polynomials that are positive on the positive orthant. Bernstein's and Handelman's theorems yield a basis to parameterize the space of polynomials that are positive on simplices and convex polytopes. Similar to SOS programming, one can find Polya's, Bernstein's and Handelman's parameterizations by solving a sequence of Linear Programs (LPs) and/or SDPs. However, unlike the SDPs associated with SOS programming, the SDPs associated with these theorems have a block-diagonal structure. This structure has been exploited in~\cite{reza_peet_tac2013} to design parallel algorithms for optimization of polynomials with large degrees and number of variables. Unfortunately, unlike SOS programming, Bernstein's, Handelman's and the original Polya's theorems do not parameterize polynomials with zeros in the positive orthant. Yet, there exist some variants of Polya's theorem which considers zeros at the corners~\cite{polya_corner} and edges~\cite{polya_edge} of simplices. Moreover, there exist other variants of Polya's theorem which provide certificates of positivity on hypercubes~\cite{peres_multisimplex,reza_CDC_hypercube}, intersection of semi-algebraic sets and the positive orthant~\cite{polya_positivstellensatz_2014} and the entire $\mathbb{R}^n$~\cite{polya_Rn}, or apply to polynomials with rational exponents~\cite{polya_rational}. 

We organize this paper as follows. In Section~\ref{sec:history}, we place Polya's, Bernstein's, Handelman's and the Positivstellensatz results in the broader topic of research on polynomial positivity. In Section~\ref{sec:poly_optim}, we first define polynomial optimization and optimization of polynomials. Then, we formulate optimization of polynomials as the problem of verifying the feasibility of semi-algebraic sets. To verify the feasibility of different semi-algebraic sets, we present algorithms based on the different variants of Polya's, Bernstein's, Handelman's and Positivstellensatz results. In Section~\ref{sec:applications}, we discuss how these algorithms apply to robust stability analysis~\cite{peres_ramos_2002,reza_peet_tac2013,peres_TAC2007} and nonlinear stability~\cite{CPA_sigurdur,reza_CDC_2013,Handelman_Sankaranarayanan,reza_CDC_2014}.  
Finally, one contribution of this paper is to apply Polya's algorithm to the problem of $H_{\infty}$ control synthesis for systems with parametric uncertainties.


\section{Background on positivity of polynomials}
\label{sec:history}	
In 1900, Hilbert published a list of mathematical problems, one of which was: For every non-negative $f \in \mathbb{R}[x]$, does there exist some non-zero $q \in \mathbb{R}[x]$ such that $q^2f$ is a sum of squares? In other words, is every non-negative polynomial a sum of squares of rational functions? This question was motivated by his earlier works~\cite{hilbert1,hilbert2}, in which he proved: 1- Every non-negative bi-variate degree 4 \textit{homogeneous} polynomial (A polynomial whose monomials all have the same degree) is a SOS of three polynomials. 2- Every bi-variate non-negative polynomial is a SOS of four rational functions. 3- Not every homogeneous polynomial with more than two variables and degree greater than 5 is SOS of polynomials. Eighty years later, Motzkin constructed a non-negative degree 6 polynomial with three variables which is not SOS~\cite{motzkin}:
\[
M(x_1,x_2,x_3)=x_1^4x_2^2+x_1^2x_2^4-3x_1^2x_2^2x_3^2+x_3^6.
\] 
Robinson~\cite{Hil17_reznick} generalized Motzkin's example as follows.  Polynomials of the form $(\prod_{i=1}^n x_i^2) f(x_1,\cdots,x_n)+1$ are not SOS if polynomial $f$ of degree $<2n$ is not SOS. Hence, although the non-homogeneous Motzkin polynomial $M(x_1,x_2,1)=x_1^2 x_2^2(x_1^2+x_2^2-3)+1$ is non-negative it is not SOS.

In 1927, Artin answered Hilbert's problem in the following theorem~\cite{artin}.
\begin{theorem} (Artin's theorem)
A polynomial $f\in \mathbb{R}[x]$ satisfies $f(x) \geq 0$ on $\mathbb{R}^n$ if and only if there exist SOS polynomials $N$ and $D \neq 0$ such that $f(x)= \frac{
N(x)}{D(x)}$.
\end{theorem}
Although Artin settled Hilbert's problem, his proof was neither constructive nor gave a characterization of the numerator $N$ and denominator $D$. In 1939, Habicht~\cite{habicht} showed that if $f$ is positive definite and can be expressed as $f(x_1, \cdots, x_n)$ $=g(x_1^2,\cdots,x_n^2)$ for some polynomial $g$, then one can choose the denominator $D=\sum_{i=1}^n x_i^2$. Moreover, he showed that by using $D=\sum_{i=1}^n x_i^2$, the numerator $N$ can be expressed as a sum of squares of monomials. Habicht used Polya's theorem (\cite{polya_book}, Theorem 56) to obtain the above characterizations for $N$ and $D$.
\begin{theorem}(Polya's theorem)
Suppose a homogeneous polynomial $p$ satisfies $p(x) > 0$ for all $x \in \{ x \in \mathbb{R}^n : x_i \geq 0, \sum_{i=1}^n {x_i} \neq 0 \}$. Then $p(x)$ can be expressed as
\[
p(x) = \dfrac{N(x)}{D(x)},
\]
where $N(x)$ and $D(x)$ are homogeneous polynomials with all positive coefficients. For every homogeneous $p(x)$ and some $e \geq 0$, the denominator $D(x)$ can be chosen as $(x_1+ \cdots +x_n)^e$.
\label{thm:polya}
\end{theorem}
Suppose $f$ is homogeneous and positive on the positive orthant and can be expressed as $f(x_1, \cdots, x_n)=g(x_1^2,\cdots,x_n^2)$ for some homogeneous polynomial $g$. By using Polya's theorem $g(y)=\frac{N(y)}{D(y)}$, where $y:=(y_1,\cdots,y_n)$ and polynomials $N$ and $D$ have all positive coefficients. By Theorem~\ref{thm:polya} we may choose $D(y)=\left( \sum_{i=1}^n y_i \right)^e$. Then $\left( \sum_{i=1}^n y_i \right)^e g(y)=N(y)$. Now let $x_i=\sqrt{y_i} $, then $\left( \sum_{i=1}^n x_i^2 \right)^e f(x_1,\cdots,x_n)=N(x_1^2,\cdots,x_n^2)$. Since $N$ has all positive coefficients, $N(x_1^2,\cdots,x_n^2)$ is a sum of squares of monomials.
Unlike the case of positive definite polynomials, it is shown that there exists no single SOS polynomial $D \neq 0$ which satisfies $f=\frac{N}{D}$ for every positive semi-definite $f$ and some SOS polynomial $N$~\cite{reznick_no_denominator}.

As in the case of positivity on $\mathbb{R}^n$, there has been an extensive research regarding positivity of polynomials on bounded sets. A pioneering result on local positivity is Bernstein's theorem (1915)~\cite{bernstein_1915}. Bernstein's theorem uses the polynomials $h_{i,j}=(1+x)^i(1-x)^j$ as a basis to parameterize univariate polynomials which are positive on $[-1,1]$.
\begin{theorem}(Bernstein's theorem)
If a polynomial $f(x) > 0$ on $[-1,1]$, then there exist $c_{i,j} > 0$ such that \[
f(x)= \sum_{ \substack{i,j \in \mathbb{N} \\ i+j=d}} c_{i,j} (1+x)^i(1-x)^j
\]
for some $d > 0$.
\end{theorem}
 Reference~\cite{positive_on_interval_reznick_powers} uses Goursat's transform of $f$ to find an upper bound on $d$. The bound is a function of the minimum of $f$ on $[-1,1]$. However, computing the minimum itself is intractable.
In 1988, Handelman~\cite{handelman} used products of affine functions as a basis (the Handelman basis) to extend Bernstein's theorem to multi-variate polynomials which are positive on convex polytopes. 
\begin{theorem} (Handelman's Theorem)
\label{thm:Handelman} Given $w_i \in \mathbb{R}^n$ and $u_i \in \mathbb{R}$, define the polytope $\Gamma^K := \{ x\in \mathbb{R}^n : w_i^Tx + u_i\geq 0, i=1,\cdots,K \}$. If a polynomial $f(x) > 0$ on $\Gamma^K$, then there exist $b_\alpha \geq 0$, $\alpha \in \mathbb{N}^K$ such that for some $d \in \mathbb{N}$,
\begin{equation}
f(x) = \sum _{\substack{\alpha \in \mathbb{N}^K \\ \alpha_1+\cdots+\alpha_K \leq d}} b_\alpha (w_1^T x+u_1)^{\alpha_1} \cdots (w_K^T x+u_K)^{\alpha_K}.
\label{eq:handelman_representation}
\end{equation}
\end{theorem}

In~\cite{roy}, first the \textit{standard triangulation} of a simplex (the convex hull of vertices in $\mathbb{R}^n$) is developed to decompose an arbitrary simplex into sub-simplices. Then, an algorithm is proposed to ensure positivity of a polynomial $f$ on the simplex by finding an expression of Form~\eqref{eq:handelman_representation} for $f$  on each sub-simplex. An upper bound on the degree $d$ in~\eqref{eq:handelman_representation} was provided in~\cite{leroy} as a function of the minimum of $f$ on the simplex, the number of variables of $f$, the degree of $f$ and the maximum of certain~\cite{leroy} affine combinations of the coefficients $b_\alpha$. Reference~\cite{roy} also provides a bound on $d$ as a function of $\max_\alpha b_\alpha$ and the minimum of $f$ over the polytope.

An extension of Handelman's theorem was made by Schweighofer~\cite{schweighofer_nonnegativity} to verify non-negativity of polynomials over compact semi-algebraic sets. Schweighofer used the cone of polynomials in~\eqref{eq:cone}
to parameterize any polynomial $f$ which has the following properties:
\begin{enumerate}
\item $f$ is non-negative over the compact semi-algebraic set $S$
\item $f=q_1 p_1 + q_2 p_2 + \cdots$ for some $q_i$ in the cone~\eqref{eq:cone} and for some $p_i > 0$ over $S  \cap \{x \in \mathbb{R}^n:f(x)=0 \}$
\end{enumerate}

\begin{theorem}(Schweighofer's theorem)
Suppose 
\[S:=\{x \in \mathbb{R}^n: g_i(x) \geq 0, g_i \in \mathbb{R}[x] \text{ for } i=1, \cdots, K\}
\]
is compact. Define the following set of polynomials which are positive on $S$.
\begin{equation}
\Theta_d:= \left\lbrace  \sum_{\substack{\lambda \in \mathbb{N}^K: \lambda_1+\cdots+\lambda_K \leq d}} s_\lambda g_1^{\lambda_1} \cdots g_K^{\lambda_K}: s_\lambda \text{ are SOS }   \right\rbrace
\label{eq:cone}
\end{equation}
If $f \geq 0$ on $S$ and there exist $q_i \in \Theta_d$ and polynomials $p_i > 0$ on $S  \cap \{x \in \mathbb{R}^n:f(x)=0 \}$ such that $f= \sum_{i} q_i p_i$
for some $d$, then $f \in \Theta_d$.
\label{Schweighofer}
\end{theorem}
 On the assumption that $g_{i}$ are affine functions, $p_i=1$ and $s_\lambda$ are constant, Schweighofer's theorem gives the same parameterization of $f$ as in Handelman's theorem. Another special case of Schweighofer's theorem is when $\lambda \in \{ 0,1\}^K$. In this case, Schweighofer's theorem reduces to 
Schmudgen's Positivstellensatz~\cite{schmudgen}. Schmudgen's Positivstellensatz states that the cone
\begin{equation}
\Lambda_g:=\left\lbrace \sum_{\lambda \in \{0,1 \}^K }s_{\lambda} g_1^{\lambda_1} \cdots g_K^{\lambda_K}: s_{\lambda} \text{ are } SOS \right\rbrace \subset \Theta_d
\label{schmudgen}
\end{equation}
is sufficient to parameterize every $f > 0$ over the semi-algebraic set $S$ generated by $\{g_1,\cdots,g_K\}$. Unfortunately, the cone $\Lambda_g$ contains $2^K$ products of $g_i$, thus finding a representation of Form~\eqref{schmudgen} for $f$ requires a search for at most $2^K$ SOS polynomials.
%
Putinar's Positivstellensatz~\cite{putinar} reduces the complexity of Schmudgen's parameterization in the case where the quadratic module of $g_i$ defined in~\eqref{eq:putinar_cone} is \textit{Archimedean}, i.e., for every $p \in \mathbb{R}[x]$, there exist $N \in \mathbb{N}$ such that
\[
N \pm p \in M_g.
\]
\begin{theorem}(Putinars's Positivstellensatz)
Let $S:=\{x \in \mathbb{R}^n: g_i(x) \geq 0, g_i \in \mathbb{R}[x] \text{ for } i=1, \cdots, K\}$ and define 
\begin{equation}
M_g:=\left\lbrace s_0+\sum_{i=1}^K s_i g_i: s_i \text{ are SOS}  \right\rbrace.
\label{eq:putinar_cone}
\end{equation}
If there exist some $N > 0$ such that $N-\sum_{i=1}^n x_i^2 \in M_g$, then $M_g$ is Archimedean. If $M_g$ is Archimedean and $f > 0$ over $S$, then $f \in M_g$.
\label{putinar}
\end{theorem}
Finding a representation of Form~\eqref{eq:putinar_cone} for $f$, only requires a search for $K+1$ SOS polynomials using SOS programming. Verifying the Archimedian condition $N-\sum_{i=1}^n x_i^2 \in M_g$ in Theorem~\ref{putinar} is also a SOS program. Observe that the Archimedian condition implies the compactness of $S$. The following theorem, lifts the compactness requirement for the semi-algebraic set $S$. 

\begin{theorem}(Stengle's Positivstellensatz)
Let $S:=\{x \in \mathbb{R}^n: g_i(x) \geq 0, g_i \in \mathbb{R}[x] \text{ for } i=1, \cdots, K\}$ and define 
\[
\Lambda_g:=\left\lbrace \sum_{\lambda \in \{0,1 \}^K }s_{\lambda} g_1^{\lambda_1} \cdots g_K^{\lambda_K}: s_{\lambda} \text{ are } SOS \right\rbrace.
\]
If $f>0$ on $S$, then there exist $p,g \in \Lambda_g$ such that $qf=p+1$.
\label{stengle}
\end{theorem}
Notice that the Parameterziation~\eqref{eq:handelman_representation} in Handelman's theorem is affine in $f$ and the coefficients $b_\alpha$. Likewise, the parameterizations in Theorems~\ref{Schweighofer} and~\ref{putinar}, i.e., $f=\sum_{\lambda} s_\lambda g_1^{\lambda_1} \cdots g_K^{\lambda_K}$ and $f=s_0+\sum_i s_i g_i$ are affine in $f,s_\lambda$ and $s_i$. Thus, one can use convex optimization to find $b_\alpha$, $s_\lambda,s_i$ and $f$. Unfortunately, since the parameterization $qf=p+1$ in Stengle's Positivstellensatz is non-convex (bilinear in $q$ and $f$), it is more difficult to verify the feasibility of $qf=p+1$ compared to Handelman's and Putinar's parameterizations.

For a comprehensive discussion on the Positivstellensatz and other polynomial positivity results in algebraic geometry see~\cite{laurent_survey,scheiderer_survey,delzell_book}.


\section{Algorithms for Polynomial Optimization}
\label{sec:poly_optim}
In this Section, we first define polynomial optimization, optimization of polynomials and its equivalent feasibility problem using semi-algebraic sets. Then, we introduce some algorithms to verify the feasibility of different semi-algebraic sets. We observe that combining these algorithms with bisection yields some lower bounds on optimal objective values of polynomial optimization problems.

\subsection{Polynomial Optimization and optimization of polynomials}
\label{sec:poly_optim2}
We define \textit{polynomial optimization} problems as 
\begin{align}
& \beta^* =  \min_{x \in \mathbb{R}^n}  \;\;\, f(x) \nonumber \\
& \text{subject to } \;\; g_i(x) \geq 0 \text{ for } i=1,\cdots,m \nonumber \\
& \hspace{0.72in} h_j(x)=0 \text { for } j=1,\cdots,r,
\label{eq:polynomial_optimization}
\end{align}
where $f,g_i,h_j \in \mathbb{R}[x]$ are given. For example, the integer program 
\begin{align}
\min_{x \in \mathbb{R}^n}  \quad &p(x) \nonumber  \\
\text{subject to} \quad &a_i^T x \geq b_i \text{ for } i=1, \cdots, m, \nonumber \\
&x \in \{-1,1\}^n,
\label{eq:polynomial_optimization_exp}
\end{align}
with given $a_i \in \mathbb{R}^n, b_i \in \mathbb{R}$ and $p \in \mathbb{R}[x]$, can be formulated as a polynomial optimization problem by setting
\begin{align*}
& f=p \\
& g_i(x)= a_i^Tx-b_i && \hspace{-1.1in} \text{ for } i=1,\cdots,m \\
& h_j(x)=x_j^2-1     && \hspace{-1.1in} \text{ for } j=1,\cdots,n.
\end{align*}
We define \textit{Optimization of polynomials} problems as
\begin{align}
&\gamma^*=\max_{x \in \mathbb{R}^n} \quad c^Tx \nonumber \\
&\text{subject to } \;\;\, F(x,y):= F_0(y)+\sum_{i=1}^n x_i F_i(y) \geq 0 \text{ for all } y \in \mathbb{R}^m,
\label{eq:optimization_of_polynomials}
\end{align}
where $c \in \mathbb{R}^n$ are given and where $F_i(y) := \sum_{\alpha \in E_d} F_{i,\alpha} y_1^{\alpha_1} \cdots y_m^{\alpha_m}$ with $E_{d_i}:=\{ \alpha \in \mathbb{N}^m : \sum_{i=1}^m \alpha_i = d_i \}$, where the coefficients $F_{i,\alpha} \in \mathbb{R}^{q \times q}$ are either given or are decision variables. Optimization of polynomials can be used to find $\beta^*$ in~\eqref{eq:polynomial_optimization}. For example, we can compute the optimal objective value $\alpha^*$ of the polynomial optimization problem
\begin{align*}
\alpha^* = \min_{x \in \mathbb{R}^n} \quad &p(x)  \nonumber  \\
\text{subject to} \quad & a_i^Tx-b_i \geq 0 && \hspace{-0.9in} \text{for } i=1, \cdots, m, \nonumber \\
&x_j^2-1=0 && \hspace{-0.9in} \text{for } j=1,\cdots,n ,
\end{align*}
by solving the problem
\begin{align}
&\alpha^*= \max_{\alpha} \quad \alpha \nonumber \\
&\text{subject to} \quad p(x) \geq \alpha \text{ for all } x \in \{-1,1\}^n \nonumber \\
&\qquad \qquad \quad \;\, a_i^T x \geq b_i \text{ for } i=1, \cdots, m \text{ and for all } x \in \{-1,1\}^n,
\label{eq:optimization_of_polynomials_exp}
\end{align}
where Problem~\eqref{eq:optimization_of_polynomials_exp} can be expressed in the Form~\eqref{eq:optimization_of_polynomials} by setting
\[
c=1, \quad n=1, \quad k=0, \quad q=m+1, \quad h_j(y)=y_j^2-1 \text{ for } j=1, \cdots,n
\]
\[
F_0(y)=\begin{bmatrix}
p(y) & 0 & \cdots & 0 \\
0 & a_1^Ty-b_1 &  & \vdots \\ 
\vdots &  & \ddots & 0 \\ 
0 & \cdots & 0 & a_m^Ty-b_m
\end{bmatrix} , \;\; 
F_1=\begin{bmatrix}
-1 & 0 & \cdots & 0 \\ 
0 & 0 &  &  \\ 
\vdots &  & \ddots & \vdots \\ 
0 &  & \cdots & 0
\end{bmatrix}.
\]

Optimization of polynomials~\eqref{eq:optimization_of_polynomials} can be formulated as the following feasibility problem. 
\begin{align}
& \hspace{-0.09in} \gamma^* = \min_\gamma \;\,  \gamma \nonumber \\
& \hspace{-0.13in} \text{ subject to }  S_\gamma:= \left\lbrace   x,y \in \mathbb{R}^n : c^T x > \gamma, \,  F(x,y) \geq 0 \right\rbrace = \emptyset,
\label{eq:feasibility}
\end{align}
where $c, F,g_i$ and $h_j$ are given. The question of feasibility of a semi-algebraic set is NP-hard~\cite{blum}. However, if we have a test to verify $S_\gamma = \emptyset$, we can find $\gamma^*$ by performing a bisection on $\gamma$.
In Section~\ref{sec:algorithms}, we use the results of Section~\ref{sec:history} to provide sufficient conditions, in the form of Linear Matrix Inequalities (LMIs), for $S_\gamma = \emptyset$.


\subsection{Algorithms}
\label{sec:algorithms}
In this section, we discuss how to find lower bounds on $\beta^*$ for different classes of polynomial optimization problems. The results in this section are primarily expressed as methods for verifying $S_{\gamma}=\emptyset$ and can be used with bisection to solve polynomial optimization problems. \vspace{0.5in}


\textbf{ \\ \noindent Case 1. Optimization over the standard simplex $\Delta^n$}

Define the standard unit simplex as
\begin{equation}
\Delta^n : = \{ x \in \mathbb{R}^n: \sum_{i=1}^n x_i=1, x_i \geq 0\}.
\label{eq:simplex}
\end{equation}
Consider the polynomial optimization problem
\begin{equation*}
\gamma^* = \min_{x \in \Delta^n} \quad f(x),
\end{equation*}
where $f$ is a homogeneous polynomial of degree $d$. If $f$ is not homogeneous, we can homogenize it by multiplying each monomial $x_1^{\alpha_1} \cdots x_n^{\alpha_n}$ in $f$ by $(\sum_{i=1}^n x_i)^{d-\Vert \alpha \Vert_1}$. Notice that since $\sum_{i=1}^n x_i=1$ for all $x \in \Delta^n$, the homogenized $f$ is equal to $f$ for all $x \in \Delta^n$. To find $\gamma^*$, one can solve the following optimization of polynomials problem. 
\begin{align}
& \gamma^* = \max_{\gamma \in \mathbb{R}} \quad \gamma \nonumber  \\
&   \text{s.t. } f(x) \geq \gamma \text{ for all } x \in \Delta^n
\label{eq:optimization_of_polynomials_simplex}
\end{align}
It can be shown that Problem~\eqref{eq:optimization_of_polynomials_simplex} is equivalent to the feasibility problem 
\begin{align*}
&\gamma^*=\min_{\gamma \in \mathbb{R}} \quad \gamma \\
&\text{s.t. } S_\gamma  := \{x \in \mathbb{R}^n: f(x) - \gamma < 0, \sum_{i=1}^n x_i=1,\, x_i \geq 0 \} = \emptyset.
\end{align*}
For a given $\gamma$, we use the following version of Polya's theorem to verify $S_\gamma = \emptyset$.

\begin{theorem}(Polya's theorem, simplex version)
If a homogeneous matrix-valued polynomial $F$ satisfies $F(x) > 0$ for all $x \in \Delta^n := \{ x \in \mathbb{R}^n: \sum_{i=1}^n x_i = 1, x_i \geq 0 \}$, then there exists $e \geq 0$ such that all the coefficients of
\[
\left( \sum_{i=1}^n x_i  \right)^e F(x)
\]
are positive definite.
\label{thm:polya_simplex}
\end{theorem}

The converse of the theorem only implies $F \geq 0$ over the unit simplex.  Given $\gamma \in \mathbb{R}$, it follows from the converse of Theorem~\ref{thm:polya_simplex} that $S_\gamma=\emptyset$ if there exist $e \geq 0$ such that 
\begin{equation}
\left( \sum_{i=1}^n x_i \right)^e \left( f(x)-\gamma \left(\sum_{i=1}^n x_i\right)^d \right)
\label{eq:polya_product}
\end{equation}
has all positive coefficients, where recall that $d$ is the degree of $f$. We can compute lower bounds on $\gamma^*$ by performing a bisection on $\gamma$. For each $\gamma$ of the bisection, if there exist $e \geq 0$ such that all of the coefficients of~\eqref{eq:polya_product} are positive, then $\gamma \leq \gamma^*$. \\

%

\textbf{ \hspace*{-0.27in} Case 2. Optimization over the hypercube $\Phi^n$:}

Given $r_i \in \mathbb{R}$, define the hypercube 
\begin{equation}
\Phi^n:=\{x \in \mathbb{R}^n:  \vert x_i \vert \leq r_i, i=1, \cdots, n  \}.\
\label{eq:hypercube}
\end{equation}
Define the set of $n$-variate \textit{multi-homogeneous} polynomials of degree vector $d \in \mathbb{N}^n$ as 
\begin{equation}
\left\lbrace
p \in \mathbb{R}[x,y]:
p(x,y) = \sum_{\substack{ h,g \in \mathbb{N}^n  \\ h+g=d  }} p_{h,g} x_1^{h_1}y_1^{g_1} \cdots x_n^{h_n}y_n^{g_n},\, p_{h,g} \in \mathbb{R}
\right\rbrace.
\label{eq:multi-homog_poly}
\end{equation}
In a more general case, if the coefficients $p_{h,g}$ are matrices, we call $p$ a \textit{matrix-valued multi-homogeneous polynomial}.
It is shown in~\cite{reza_CDC_hypercube} that for every polynomial $f(z)$ with $z \in \Phi^n$, there exists a multi-homogeneous polynomial $p$ such that 
\begin{equation}
\left\lbrace f(z) \in \mathbb{R} : z \in \Phi^n \right\rbrace = \left\lbrace p(x, y) \in \mathbb{R} : \, x, y \in \mathbb{R}^n \text{ and } (x_i, y_i) \in \Delta^2 \text{ for } i=1, \cdots,n \right\rbrace. 
\label{eq:fz_multi}
\end{equation}
For example, consider $f(z_1,z_2)=z_1^2+z_2$, with $z_1 \in [-2,2]$ and $z_2 \in [-1,1]$. Let $x_1=\frac{z_1+2}{4} \in [0,1]$ and $x_2=\frac{z_2+1}{2}\in [0,1]$. Then define
\[
q(x_1,x_2) := f(4x_1-2,2x_2-1) =  16x_1^2-16x_1+2x_2+3
\]
By homogenizing $q$ we obtain the multi-homogeneous polynomial
\begin{align*}
p(x,y)= & 16x_1^2(x_2+y_2)-16x_1(x_1+y_1)(x_2+y_2)+2 x_2(x_1+y_1)^2 \\
&+3(x_1+y_1)^2(x_2+y_2),\;\; (x_1,y_1), (x_2,y_2) \in \Delta^2
\end{align*}
with degree vector $d=[2,1]$, where $d_1=2$ is the sum of exponents of $x_1$ and $y_1$ in every monomial of $p$, and $d_2=1$ is the sum of exponents of $x_2$ and $y_2$ in every monomial of $p$.
See~\cite{reza_CDC_hypercube} for an algorithm which computes the multi-homogeneous polynomial $p$ for an arbitrary $f$ defined on a hypercube.

Now consider the polynomial optimization problem
\[
\gamma^*=\min_{x \in \Phi^n} f(x).
\]
To find $\gamma^*$, one can solve the following feasibility problem.
\begin{align}
&\gamma^*=\min_{\gamma \in \mathbb{R}} \quad \gamma \nonumber \\
&\text{subject to } \;\, S_{\gamma,r}  := 
\{x \in \mathbb{R}^n: f(x) - \gamma < 0, \, \vert x_i \vert \leq r_i,\, i=1, \cdots,n \} = \emptyset
\label{eq:feasibility_hypercube}
\end{align}
For a given $\gamma$, one can use the following version of Polya's theorem to verify $S_{\gamma,r} = \emptyset$.

\begin{theorem}(Polya's theorem, multi-homogeneous version)
A matrix-valued multi-homogeneous polynomial $F$ satisfies $F(x,y) > 0$ for all $(x_i,y_i) \in \Delta^2, i=1, \cdots, n$, if there exist $e \geq 0$ such that all the coefficients of
\[
\left( \prod_{i=1}^n \left( x_i+y_i \right)^e \right) F(x,y)
\]
are positive definite.
\label{thm:polya_multi-simplex}
\end{theorem}

The Converse of the theorem only implies non-negativity of $F$ over the multi-simplex. To find lower bounds on $\gamma$, we first obtain the multi-homogeneous form $p$ of the polynomial $f$ in~\eqref{eq:feasibility_hypercube} by using the algorithm in Section 3 of~\cite{reza_CDC_hypercube} . Given $\gamma \in \mathbb{R}$ and $r \in \mathbb{R}^n$, from the converse of Theorem~\ref{thm:polya_multi-simplex} it follows that $S_{\gamma,r}=\emptyset$ in~\eqref{eq:feasibility_hypercube} if there exist $e \geq 0$ such that 
\begin{equation}
\left( \prod_{i=1}^n \left( x_i+y_i \right)^e \right) \left( p(x,y) - \gamma \left( \prod_{i=1}^n \left( x_i+y_i \right)^{d_i} \right) \right)
\label{eq:polya_product_hypercube}
\end{equation}
has all positive coefficients, where $d_i$ is the degree of $x_i$ in $p(x,y)$. We can compute lower bounds on $\gamma^*$ by performing a bisection on $\gamma$. For each $\gamma$ of the bisection, if there exist $e \geq 0$ such that all of the coefficients of~\eqref{eq:polya_product_hypercube} are positive, then $\gamma \leq \gamma^*$.\\


%

\textbf{ \hspace*{-0.27in} Case 3. Optimization over the convex polytope $\Gamma^K$:}

Given $w_i \in \mathbb{R}^n$ and $u_i \in \mathbb{R}$, define the convex polytope 
\[
\Gamma^K := \{ x\in \mathbb{R}^n : w_i^Tx + u_i \geq 0, i=1,\cdots,K \}.
\]
 Suppose $\Gamma^K$ is bounded. Consider the polynomial optimization problem
\[
\gamma^*=\min_{x \in \Gamma^K} f(x),
\]
where $f$ is a polynomial of degree $d_f$. To find $\gamma^*$, one can solve the feasibility problem.
\begin{align*}
&\gamma^*=\min_{\gamma \in \mathbb{R}} \quad \gamma \nonumber \\
&\text{subject to } \;\, S_{\gamma,K}  := 
\{x \in \mathbb{R}^n: f(x) - \gamma <  0, \, w_i^Tx+u_i \geq 0, \, i = 1, \cdots, K \} = \emptyset.
\end{align*}
Given $\gamma$, one can use Handelman's theorem (Theorem~\ref{thm:Handelman}) to verify $S_{\gamma,K}=\emptyset$ as follows. 
Consider the Handelman basis associated with polytope $\Gamma^K$ defined as
\begin{equation*}
B_s:=\left\lbrace \lambda_{\alpha} \in \mathbb{R}[x] : \lambda_{\alpha}(x)= \prod_{i=1}^K \left( w_i^Tx+u_i \right)^{\alpha_i}, \alpha \in \mathbb{N}^K, \sum_{i=1}^K  \alpha_i \leq s \right\rbrace.
\end{equation*} 
Basis $B_s$ spans the space of polynomials of degree $s$ or less, however it is not minimal. Given polynomial $f(x)$ of degree $d_f, \gamma \in \mathbb{R}$ and $d_{\max} \in \mathbb{N}$, if there exist 
\begin{equation}
c_\alpha \geq 0 \text{ for all } \alpha\in \{\alpha \in \mathbb{N}^K: \Vert \alpha \Vert_1 \leq d \}
\label{eq:handelman_representation1}
\end{equation}
such that
\begin{equation}
f(x) - \gamma = \sum_{\Vert \alpha \Vert_1 \leq d} c_\alpha \prod_{i=1}^K (w_i^Tx+u_i)^{\alpha_i},
\label{eq:handelman_representation2}
\end{equation}
for $d=d_f$, then $f(x) - \gamma \geq 0$ for all $x \in \Gamma^K$. Thus $S_{\gamma,K} = \emptyset$. Feasibility of Conditions~\eqref{eq:handelman_representation1} and~\eqref{eq:handelman_representation2} can be determined using linear programming. If~\eqref{eq:handelman_representation1} and~\eqref{eq:handelman_representation2} are infeasible for some $d$, then one can increase $d$ up to $d_{\max}$. From Handelman's theorem, if $f(x) - \gamma > 0$ for all $x \in \Gamma^K$, then for some $d \geq d_f$, Conditions~\eqref{eq:handelman_representation1} and~\eqref{eq:handelman_representation2} hold. However, computing upper bounds for $d$ is difficult~\cite{reznick_powers_polyhedra,leroy}.

Similar to Cases 1 and 2, we can compute lower bounds on $\gamma^*$ by performing a bisection on $\gamma$. For each $\gamma$ of the bisection, if there exist $d \geq d_f$ such that Conditions~\eqref{eq:handelman_representation1} and~\eqref{eq:handelman_representation2}, then $\gamma \leq \gamma^*$.\\

%

\textbf{ \hspace*{-0.27in} Case 4: Optimization over compact semi-algebraic sets:} 

Recall that we defined a semi-algebraic set as
\begin{equation}
S:=\{ x \in \mathbb{R}^n : g_i(x) \geq 0,i=1,\cdots,m, \, h_j(x) = 0,j=1,\cdots,r \}.
\label{eq:semialgebraic_set2}
\end{equation}
Suppose $S$ is compact. Consider the polynomial optimization problem
\begin{align*}
& \gamma^*= \min_{x \in \mathbb{R}^n}  \quad f(x) \nonumber \\
& \text{subject to } \;\;\, g_i(x) \geq 0 \text{ for } i=1,\cdots,m \nonumber \\
& \hspace{0.75in} h_j(x)=0 \text { for } j=1,\cdots,r.
\end{align*}
Define the following cone of polynomials which are positive over $S$.
\begin{equation}
M_{g,h} \hspace{-0.03in} := \hspace{-0.03in} \left\lbrace \hspace{-0.03in} m \hspace{-0.03in} \in \hspace{-0.03in} \mathbb{R}[x] :  m(x) \hspace{-0.03in} - \hspace{-0.04in} \sum_{i=1}^m s_i(x) g_i(x) \hspace{-0.03in} - \hspace{-0.04in} \sum_{i=1}^r t_i(x) h_i(x) \text{ is SOS}, s_i \in \Sigma_{2d} , t_i \in \mathbb{R}[x]  \right\rbrace,
\label{eq:putinar_cone2}
\end{equation}
where $\Sigma_{2d}$ denotes the set of SOS polynomials of degree $2d$. 
From Putinar's Positivstellensatz (Theorem~\ref{putinar}) it follows that if the Cone~\eqref{eq:putinar_cone2} 
is Archimedean, then the solution to the following SOS program is a lower bound on $\gamma^*$. Given $d \in \mathbb{N}$, define
\begin{align}
& \gamma^d =  \max_{\gamma \in \mathbb{R}, s_i, t_i} \;\; \gamma \nonumber \\
&\text{subject to } \; f(x)-\gamma - \sum_{i=1}^m s_i(x) g_i(x) -  \sum_{i=1}^r t_i(x) h_i(x) \text{ is SOS },\, t_i \in \mathbb{R}[x], s_i \in \Sigma_{2d}.
\label{eq:SOS_putinar}
\end{align}
For given $\gamma \in \mathbb{R}$ and $d \in \mathbb{N}$, Problem~\eqref{eq:SOS_putinar} is the following linear matrix inequality.
\begin{align}
&\text{Find } \quad\, Q_i \geq 0, \, P_j \;\; \text{ for } i=0,\cdots,m \text{ and } j=1, \cdots, r \nonumber \\ 
& \hspace*{-0.15in} \text{such that } \; f(x)-\gamma = z_d^T(x)\left( Q_0+\sum_{i=1}^m Q_i g_i(x) + \sum_{j=1}^r P_j h_j(x) \right)z_d(x),
\label{eq:putinar_feasibility}
\end{align}
where $Q_i,P_j \in \mathbb{S}^{N}$, where $\mathbb{S}^N$ is the subspace of symmetric matrices in $\mathbb{R}^{N \times N}$ and $N:=\dbinom{n+d}{d}$, and where $z_d(x)$ is the vector of monomial basis of degree $d$ or less. See~\cite{sostools2013,SDP_SIAMbook_parrilo} for methods of solving SOS programs. It is shown in~\cite{lasserre2001} that if the Cone~\eqref{eq:putinar_cone2} is Archimedean, then $ \lim_{d \to \infty} \gamma^d = \gamma^*$. 

%

If the Cone~\eqref{eq:putinar_cone2} is not Archimedean, then we can use Schmudgen's Positivstellensatz to obtain the following SOS program with solution $\gamma^d \leq \gamma^*$.
\begin{align}
& \gamma^d =  \max_{\gamma \in \mathbb{R}, s_i, t_i} \quad \gamma \nonumber \\
&\text{subject to } \; f(x)-\gamma = 1 + \hspace*{-0.15in} \sum_{\lambda \in \{0,1\}^m} \hspace*{-0.085in} s_\lambda(x) g_1(x)^{\lambda_1} \cdots g_m(x)^{\lambda_m} + \sum_{i=1}^r t_i(x) h_i(x), \, t_i \in \mathbb{R}[x],  \nonumber \\ 
& \hspace{4in} s_\lambda \in \Sigma_{2d}. 
\label{eq:SOS_schmudgen}
\end{align}

The Positivstellensatz and SOS programming can also be applied to polynomial optimization over a more general form of semi-algebraic sets defined as
\begin{align*}
&T:= \\
&\{ x\in \mathbb{R}^n: g_i(x) \geq 0,\, i =1, \cdots,m, \, h_j(x) = 0, \, j=1, \cdots,r, \, q_k(x) \neq 0, \, k=1, \cdots, l \}.
\end{align*}
It can be shown that $T = \emptyset$ if and only if 
\begin{align*}
 \hat{T}:= \{ (x,y) \in \mathbb{R}^{n+l}: g_i(x) \geq 0,\, i =1, \cdots,m, \,  h_j(x) & =  0, \, j=1, \cdots,r, \\
& y_k q_k(x) = 1, \, k=1, \cdots, l \} = \emptyset.
\end{align*}
Thus, for any $f \in \mathbb{R}[x]$, we have 
\[
\min\limits_{x\in T} f(x)= \min\limits_{(x,y) \in \hat{T}} f(x).
\]
Therefore, to find lower bounds on $\min_{x\in T} f(x)$, one can apply SOS programming and Putinar's Positivstellensatzs to $\min_{(x,y) \in \hat{T}} f(x)$. We have already addressed this problem in the current case.\\

\textbf{ \hspace*{-0.27in} Case 5: Tests for non-negativity on $\mathbb{R}^n$:} 

The following theorem~\cite{habicht}, gives a test for non-negativity of a class of homogeneous polynomials.
\begin{theorem}(Habicht theorem)
For every homogeneous polynomial $f$ that satisfies $f(x_1,\cdots,x_n) > 0$ for all $x \in \mathbb{R}^n \setminus \{0\}$ and $f(x_1,\cdots,x_n)=g(x_1^2,\cdots,x_n^2)$ for some polynomial $g$, there exist $e \geq 0$ such that all of the coefficients of
\begin{equation}
\left( \sum_{i=1}^n x_i^2 \right)^e f(x_1,\cdots,x_n)
\label{eq:habicht_product}
\end{equation}
are positive. In particular, the product is a sum of squares of monomials. 
\label{thm:habicht}
\end{theorem}
Habicht's theorem defines a test for non-negativity of any homogeneous polynomial $f$ of the form $f(x_1,\cdots,x_n)=g(x_1^2,\cdots,x_n^2)$ as follows. Multiply $f$ repeatedly by $\sum_{i=1}^n x_i^2$. If for some $e \in \mathbb{N}$, the Product~\eqref{eq:habicht_product} has all positive coefficients, then $f \geq 0$.
An alternative test for non-negativity on $\mathbb{R}^n$ is given in the following theorem~\cite{polya_Rn}. 
\begin{theorem}
Define $E_n:=\{-1,1\}^n$. Suppose a polynomial $f(x_1,\cdots,x_n)$ of degree $d$ satisfies $f(x_1,\cdots,x_n) > 0$ for all $x \in \mathbb{R}^n$ and its homogenization is positive definite. Then
\begin{enumerate}
\item there exist $\lambda_e \geq 0$ and coefficients $c_\alpha \in \mathbb{R}$ such that
\begin{equation}
\left(1+ e^Tx \right)^{\lambda_e} f(x_1, \cdots,x_n) = \sum_{\alpha \in  I_e} c_\alpha x_1^{\alpha_1} \cdots x_n^{\alpha_n} \text{ for all } e \in E_n,
\label{eq:santos_product}
\end{equation}  
where $I_e := \{ \alpha \in \mathbb{N}^n: \Vert \alpha \Vert_1 \leq d+\lambda_e \}$ and $sgn(c_\alpha) = e_1^{\alpha_1} \cdots e_n^{\alpha_n}$.

\item there exist positive $N,D \in \mathbb{R}[x_1^2,\cdots,x_n^2,f^2]$ such that $f=\frac{N}{D}$.
\end{enumerate}
\label{thm:polya_satos}
\end{theorem}
Based on the converse of Theorem~\ref{thm:polya_satos}, we can propose the following test for non-negativity of polynomials over the cone $\Lambda_e:=\{x \in \mathbb{R}^n: sgn(x_i)=e_i, i=1,\cdots,n \}$ for some $e \in E_n$. Multiply a given polynomial $f$ repeatedly by $1+e^Tx$ for some $ e \in E_n$. If there exists $\lambda_e \geq 0$ such that $sgn(c_\alpha) = e_1^{\alpha_1} \cdots e_n^{\alpha_n}$, then $f(x) \geq 0$ for all $x \in \Lambda_e$.
Since $\mathbb{R}^n = \cup_{e \in E_n} \Lambda_e$, we can repeat the test $2^n$ times to obtain a test for non-negativity of $f$ over $\mathbb{R}^n$.


The second part of Theorem~\ref{thm:polya_satos} gives a solution to the Hilbert's problem in Section~\ref{sec:history}. See~\cite{polya_Rn} for an algorithm which computes polynomials $N$ and $D$.


\section{Applications of polynomial optimization}
\label{sec:applications}
In this section, we discuss how the algorithms in Section~\ref{sec:algorithms} apply to stability analysis and control of dynamical systems. We consider robust stability analysis of linear systems with parametric uncertainty, stability of nonlinear systems, robust controller synthesis for systems with parametric uncertainty and stability of systems with time-delay.

\subsection{Robust stability analysis}
\label{sec:applications_lin}
Consider the linear system
\begin{equation}
\dot{x}(t)=A(\alpha)x(t),
\label{eq:sys_lin}
\end{equation}
where $A(\alpha) \in \mathbb{R}^{n \times n}$ is a polynomial and $\alpha \in Q \subset \mathbb{R}^l$ is the vector of uncertain parameters, where $Q$ is compact. From converse Lyapunov theory~\cite{khalil} and existence of polynomial solutions for feasible patameter-dependent LMIs ~\cite{bliman_poly_solution_2004} it follows that System~\eqref{eq:sys_lin} is asymptotically stable if and only if there exist matrix-valued polynomial $P(\alpha) \in \mathbb{S}^n$ such that
\begin{equation}
P(\alpha) > 0 \text{ and } A^T(\alpha) P(\alpha) + P(\alpha) A(\alpha) < 0 \text{ for all } \alpha \in Q.
\label{eq:LyapIneq_linear}
\end{equation}
If $Q$ is a semi-algebraic set, then asymptotic stability of System~\eqref{eq:sys_lin} is equivalent to positivity of $\gamma^*$ in the following optimization of polynomials problem for some $d \in \mathbb{N}$.
\begin{align}
&\gamma^* = \max_{\gamma \in \mathbb{R}, C_\beta \in \mathbb{S}^n} \;\; \gamma \nonumber  \\
&\text{{\small subject to}} 
\begin{small}
\begin{bmatrix}
\sum_{\beta \in E_d} C_\beta \alpha^\beta & \hspace*{-0.3in}    0 \\ 
0                                         &  \hspace*{-0.3in}  -A^T(\alpha) \hspace*{-0.04in} \left( \sum_{\beta \in E_d} C_\beta \alpha^\beta \right) \hspace*{-0.04in} - \hspace*{-0.04in} \left (\sum_{\beta \in E_d} C_\beta \alpha^{\beta} \right) \hspace*{-0.04in} A(\alpha) 
\end{bmatrix}
 \hspace*{-0.05in} - \hspace*{-0.03in}
\gamma I \alpha^T \hspace*{-0.03in} \alpha  \geq 0, \alpha \in Q,
 \end{small}
\label{eq:OOP_linear}
\end{align}
where we have denoted $\alpha_1^{\beta_1} \cdots \alpha_l^{\beta_l}$ by $\alpha^\beta$ and 
\begin{equation}
E_d:= \{ \beta \in \mathbb{N}^l: \sum_{i=1}^n \beta_i \leq d\}.
\label{eq:Ed}
\end{equation}
 Given stable systems of the form defined in~\eqref{eq:sys_lin} with different classes of polynomials $A(\alpha)$, we discuss different algorithms for solving~\eqref{eq:OOP_linear}. Solutions to~\eqref{eq:OOP_linear} yield Lyapunov functions of the form $V=x^T(\sum_{\beta \in E_d} C_\beta \alpha^\beta)x$ proving stability of System~\eqref{eq:sys_lin}.

\noindent \textbf{\\Case 1. $A(\alpha)$ is affine with $\alpha \in \Delta^l$:}\\
Consider the case where $A(\alpha)$ belongs to the polytope
\[
\Lambda_l:= \left\lbrace   A(\alpha) \in \mathbb{R}^{n \times n}: A(\alpha)=\sum_{i=1}^{l} A_i \alpha_i, A_i \in \mathbb{R}^{n \times n}, \alpha_i \in \Delta^l \right\rbrace ,
\]
where $A_i$ are the vertices of the polytope and $\Delta^l$ is the standard unit simplex defined as in~\eqref{eq:simplex}. Given $A(\alpha) \in \Lambda_l$, we address the problem of stability analysis of System~\eqref{eq:sys_lin} for all $\alpha \in \Delta^l$.

A sufficient condition for asymptotic stability of System~\eqref{eq:sys_lin} is to find a matrix $P > 0$ such that the Lyapunov inequality $A^T(\alpha) P + P A(\alpha) < 0$ holds for all $\alpha \in \Delta^l$. If $A(\alpha) = \sum_{i=1}^l A_i \alpha_i$, then from convexity of $A$ it follows that the condition
\[
A^T(\alpha) P + P A(\alpha) < 0 \text{ for all } \alpha \in \Delta^l
\]
is equivalent to positivity of $\gamma^*$ in the following semi-definite program.
\begin{align}
&\gamma^* = \max_{\gamma \in \mathbb{R},P \in \mathbb{S}^n} \;\; \gamma \nonumber \\
&\text{ subject to } 
\begin{bmatrix}
P      &  0               &  \cdots  &      0 \\ 
0      & -A_1^T P - P A_1 &  0       & \vdots \\ 
\vdots &    0             &  \ddots  &    0   \\ 
0      &      \cdots      &   0      & -A_l^T P - P A_l 
\end{bmatrix}
-
\gamma I \geq 0,
\label{eq:SDP1}
\end{align}
Any $P \in \mathbb{S}^n$ that satisfies the LMI in~\eqref{eq:SDP1} for some $\gamma > 0$, yields a Lyapunov function of the form $V=x^TPx$. However for many systems, this class of Lyapunov functions can be conservative (see Numerical Example 1). 

More general classes of Lyapunov functions such as parameter-dependent functions of the forms $V=x^T(\sum_{i=1}^l P_i \alpha_i)x$~\cite{peres_ramos_2002,gahinet_chilali_TAC1996,oliveira_peres_2005} and $V=x^T(\sum_\beta {P_{\beta \in E_d} \alpha^\beta})x$~\cite{peres_TAC2007,reza_peet_tac2013} have been utilized in the literature. 
As shown in~\cite{peres_ramos_2002}, given $A_i \in \mathbb{R}^{n \times n}$, $x^TP(\alpha)x$ with $P(\alpha) = \sum_{i=1}^l P_i \alpha_i$ is a Lyapunov function for~\eqref{eq:sys_lin} with $\alpha \in \Delta^l$ if the following LMI consitions hold.
\begin{align*}
& P_i > 0 \;\;    && \text{for } i=1, \cdots,l \\
& A_i^T P_i + P_i A_i < 0 \;\; && \text{for } i=1, \cdots,l \\
& A_i^T P_j + A_j^T P_i + P_jA_i + P_iA_j < 0 \;\; && \text{for } i=1, \cdots,l-1, \, j=i+1, \cdots, l
\end{align*}

In~\cite{bliman_homogeneous}, it is shown that given continuous functions $A_i,B_i : \Delta^l \rightarrow \mathbb{R}^{n \times n}$ and continuous function $R: \Delta^l \rightarrow \mathbb{S}^n$, if there exists a continuous function $X:\Delta^l \rightarrow \mathbb{S}^n$ which satisfies
\begin{equation}
\sum_{i=1}^N \left( A_i(\alpha) X(\alpha) B_i(\alpha) + B_i(\alpha)^T X(\alpha) A_i(\alpha)^T \right) + R(\alpha) > 0 \; \text{ for all } \alpha \in \Delta^l,
\label{eq:PD_LMI}
\end{equation}
then there exists a homogeneous polynomial $Y:\Delta^l \rightarrow \mathbb{S}^n$ which also satisfies~\eqref{eq:PD_LMI}.
Motivated by this result,~\cite{peres_TAC2007} uses the class of homogeneous polynomials of the form 
\begin{equation}
P(\alpha)= \sum_{\beta \in I_d} P_\beta \alpha_1^{\beta_1} \cdots \alpha_l^{\beta_l},
\label{eq:P(alpha)}
\end{equation}
with
\begin{equation}
I_d:= \left\lbrace \beta \in \mathbb{N}^{l}: \sum_{i=1}^l \beta_i = d \right\rbrace 
\label{eq:index_set}
\end{equation}
to provide the following necessary and sufficient LMI condition for stability of System~\eqref{eq:sys_lin}. Given $A(\alpha) = \sum_{i=1}^l A_i \alpha_i$, System~\eqref{eq:sys_lin} is asymptotically stable for all $\alpha \in \Delta^l$ if and only if there exist some $d \geq 0$ and positive definite $P_\beta \in \mathbb{S}^n, \beta \in I_d$ such that
\begin{equation}
\sum_{\substack{ i=1,\cdots,l \\ \beta_i > 0 }} \left( A_i^T P_{\beta-e_i} + P_{\beta-e_i} A_i \right) < 0  \text{ for all } \beta \in I_{d+1},
\label{eq:LMI_peres}
\end{equation}
where $e_i = [ \, 0 \cdots \; 0 \; \underbrace{1}_{i^{th}} \; 0 \cdots 0 \, ] \in \mathbb{N}^l, i=1, \cdots, l$ form the canonical basis for $\mathbb{R}^l$.

\noindent \textbf{\\Numerical Example 1:} Consider the system $\dot{x}(t) = A(\alpha,\eta) x(t)$ from~\cite{chesi_ppdlf_2005}, where $A(\alpha,\eta)= (A_0+A_1 \eta) \alpha_1 + (A_0+A_2 \eta) \alpha_2 + (A_0+A_3 \eta) \alpha_3$, where
\[
A_0 = \begin{bmatrix}
-2.4 & -0.6 & -1.7 & 3.1  \\
0.7 & -2.1  & -2.6 & -3.6 \\
0.5 &  2.4  & -5.0 & -1.6 \\
-0.6 & 2.9  & -2.0 & -0.6
\end{bmatrix} ,
A_1 = \begin{bmatrix}
1.1 & -0.6 & -0.3 & -0.1 \\
-0.8 &  0.2 & -1.1 &  2.8 \\
-1.9 &  0.8 & -1.1 &  2.0 \\
-2.4 & -3.1 & -3.7 & -0.1
\end{bmatrix} ,
\]
\[
A_2 = \begin{bmatrix}
0.9 &  3.4 &  1.7 &  1.5 \\
-3.4 & -1.4 &  1.3 &  1.4 \\
1.1 &  2.0 & -1.5 & -3.4 \\
-0.4 &  0.5 &  2.3 &  1.5 \\
\end{bmatrix} ,
A_3 = \begin{bmatrix}
-1.0 & -1.4 & -0.7 & -0.7 \\
2.1 &  0.6 & -0.1 & -2.1 \\
0.4 & -1.4 &  1.3 &  0.7 \\
1.5 &  0.9 &  0.4 & -0.5
\end{bmatrix} 
\]
and $(\alpha_1, \alpha_2, \alpha_3) \in \Delta^3, \, \eta \geq 0$. We would like to find $\eta^* = \max \eta$  such that $\dot{x}(t) = A(\alpha,\eta) x(t)$ is asymptotically stable for all $\eta \in [0,\eta^*]$.

\begin{figure}[htbp]
\includegraphics[scale=0.25]{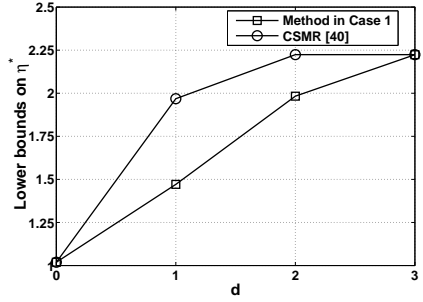}
\caption{lower-bounds for $\eta^*$ computed using the LMIs in~\eqref{eq:LMI_peres} and the method in~\cite{chesi_ppdlf_2005}}
\label{fig:simplex_e0}
\end{figure}

By performing a bisection on $\eta$ and verifying the inequalities in~\eqref{eq:LMI_peres} for each $\eta$ of the bisection algorithm, we obtained lower bounds on $\eta^*$ (see  Figure~\ref{fig:simplex_e0}) using $d=0,1,2$ and $3$. For comparison, we have also plotted the lower bounds computed in~\cite{chesi_ppdlf_2005} using the Complete Square Matricial Representation (CSMR) of the Lyapunov inequalities in~\eqref{eq:LyapIneq_linear}. Both methods found $\max \eta=2.224$, however the method in~\cite{chesi_ppdlf_2005} used a lower $d$ to find this bound.

\noindent \textbf{\\Case 2. $A(\alpha)$ is a polynomial with $\alpha \in \Delta^l$:}

Given $A_h \in \mathbb{R}^{n \times n}$ for $h \in I_d$ as defined in~\eqref{eq:index_set}, we address the problem of stability analysis of System~\eqref{eq:sys_lin} with $A(\alpha)=\sum_{h \in I_d} A_h \alpha_1^{h_1} \cdots \alpha_l^{h_l}$ for all $\alpha \in \Delta^l$.
 Using Lyapunov theory, this problem can be formulated as the following optimization of polynomials problem.
\begin{align}
&\gamma^*=\max_{\gamma \in \mathbb{R}, P \in \mathbb{R}[\alpha]} \;\; \gamma \nonumber \\
& \text{subject to } \begin{bmatrix}
P(\alpha)      &  0           \\ 
0      & -A(\alpha)^T P(\alpha) - P(\alpha) A(\alpha) 
\label{eq:OOP_simplex} 
\end{bmatrix}
-
\gamma I \geq 0  \text{ for all } \alpha \in \Delta^l
\end{align}
System~\eqref{eq:sys_lin} is asymptotically stable for all $\alpha \in \Delta^l$ if and only if $\gamma^* > 0$. As in Case 1 of Section~\ref{sec:algorithms}, one can apply bisection algorithm on $\gamma$ and use Polya's theorem (Theorem~\ref{thm:polya_simplex}) as a test for feasibility of Constraint~\eqref{eq:OOP_simplex} to find lower bounds on $\gamma^*$.
Suppose $P$ and $A$ are homogeneous matrix valued polynomials. Given $\gamma \in \mathbb{R}$, it follows from Theorem~\ref{thm:polya_simplex} that the inequality condition in~\eqref{eq:OOP_simplex} holds for all $\alpha \in \Delta^l$ if there exist some $e \geq 0$ such that
\begin{align}
  &\left( \sum_{i=1}^l \alpha_i \right)^e \left( P(\alpha) - \gamma I \left( \sum_{i=1}^l \alpha_i \right)^{d_p} \right) 
\label{eq:polya_ineq0} \\
  & \hspace{-0.8in} \text{ and } \nonumber \\
- &\left( \sum_{i=1}^l \alpha_i \right)^e  \left( A(\alpha)^T  P(\alpha) 
 +  P(\alpha)  A(\alpha) + \gamma I \left( \sum_{i=1}^l \alpha_i \right)^{d_p+d_a} \right) 
 \label{eq:polya_ineq1}
\end{align} 
have all positive coefficients, where $d_p$ is the degree of $P$ and $d_a$ is the degree of $A$. Let $P$ and $A$ be of the form
\begin{equation}
P(\alpha)= \sum_{h \in I_{d_p}} P_h \alpha_1^{h_1} \cdots \alpha_l^{h_l}, P_h \in \mathbb{S}^n \;\; \text{ and } \;\;
A(\alpha) = \sum_{h \in I_{d_a}} A_{h} \alpha_1^{h_1} \cdots \alpha_l^{h_l}, A_h \in \mathbb{R}^{n \times n}.
\label{eq:P_A_simplex}
\end{equation}
By combining~\eqref{eq:P_A_simplex} with~\eqref{eq:polya_ineq0} and~\eqref{eq:polya_ineq1} it follows that for a given $\gamma \in \mathbb{R}$, the inequality condition in~\eqref{eq:OOP_simplex} holds for all $\alpha \in \Delta^l$ if there exist some $e \geq 0$ such that
\begin{align}
 & \hspace*{-0.18in} \left( \sum_{i=1}^l \alpha_i \right)^e \hspace{-0.02in} \left( \sum_{h \in I_{d_p}} \hspace*{-0.03in} P_h \alpha_1^{h_1} \cdots \alpha_l^{h_l}  - \gamma I \left( \sum_{i=1}^l \alpha_i \right)^{d_p} \right) = \hspace{-0.1in}  \sum_{g \in I_{d_p+e}} \hspace*{-0.05in} \left( \sum_{h \in I_{d_p}} f_{g,h} P_h \right) \alpha_1^{g_1} \cdots \alpha_l^{g_l} \label{eq:prod_polya_simplex1} \\
 & \hspace*{-0.16in}  \text{  and } \nonumber \\
 - &\left( \sum_{i=1}^l \alpha_i \right)^e \hspace*{-0.05in} \left( \hspace*{-0.06in} \left( \sum_{h \in I_{d_a}} \hspace*{-0.05in} A^T_{h} \alpha^h \right) \hspace*{-0.06in} \left( \sum_{h \in I_{d_p}} \hspace*{-0.05in} P_h \alpha_1^{h_1} \cdots \alpha_l^{h_l} \right) \hspace*{-0.05in} + \hspace*{-0.05in} \left( \sum_{h \in I_{d_p}} \hspace*{-0.05in} P_h \alpha_1^{h_1} \cdots \alpha_l^{h_l} \hspace*{-0.04in} \right) \hspace*{-0.07in} \left( \sum_{h \in I_{d_a}} \hspace*{-0.05in} A_{h} \alpha^h\hspace*{-0.04in} \right)  \right. \nonumber \\
& \hspace*{-0.18in} \left.   + \gamma I \left( \sum_{i=1}^l \alpha_i \right)^{d_p+d_a} \hspace*{-0.02in} \right) = \hspace*{-0.09in} \sum_{q \in I_{d_a+d_p+e}} \hspace*{-0.07in} \left( \sum_{h \in I_{d_p}} M_{h,q}^T P_h + P_h M_{h,q}  \hspace*{-0.04in} \right)\alpha_1^{q_1} \cdots \alpha_l^{q_l}
\label{eq:prod_polya_simplex2}
\end{align}
have all positive coefficients, i.e.,
\begin{align}
& \sum_{h \in I_{d_p}} f_{h,g} P_h > 0 \;\; \text{ for all } g \in I_{d_p+e} \nonumber \\
& \sum_{h \in I_{d_p}} \left(  M^T_{h,q} P_h + P_h M_{h,q}   \right) < 0 \;\; \text{ for all } q \in I_{d_p+d_a+e},
\label{eq:polya_simplex_conditions}
\end{align}
where we define $f_{h,g} \in \mathbb{R}$ as the coefficient of $P_h \alpha_1^{g_1} \cdots \alpha_l^{g_l}$ after expanding~\eqref{eq:prod_polya_simplex1}. Likewise, we define $M_{h,q} \in \mathbb{R}^{n \times n}$ as the coefficient of $P_h \alpha_1^{q_1} \cdots \alpha_l^{q_l}$ after expanding~\eqref{eq:prod_polya_simplex2}. See~\cite{reza_peet_acc_2012} for recursive formulae for $f_{h,g}$ and $M_{h,q}$. Feasibility of Conditions~\eqref{eq:polya_simplex_conditions} can be verified by the following semi-definite program.
\begin{small}
\begin{align}
& \max_{\eta \in \mathbb{R}, P_h \in \mathbb{S}_+^n} \eta \nonumber \\
& \text{subject to} \begin{bmatrix}
\sum_{h \in I_{d_p}} \hspace*{-0.03in} f_{h,g^{(1)}} P_h & 0 &  & \ldots &  & 0 \\ 
0 & \ddots &  &  &  &  \\ 
 &  &  & \hspace*{-0.6in} \sum_{h \in I_{d_p}} f_{h,g^{(L)}} P_h &  & \vdots \\ 
\vdots &  &  &  & \hspace*{-0.8in} -\sum_{h \in I_{d_p}} \left(M^T_{h,q^{(1)}} P_h + P_h M_{h,q^{(1)}} \right) &  \\ 
 &  &  &  & \hspace*{-0.5in} \ddots & 0 \\ 
0 & \cdots & 0  &  &  & \hspace*{-1.74in} -\sum_{h \in I_{d_p}} \left(M^T_{h,q^{(M)}} P_h + P_h M_{h,q^{(M)}} \right)  \hspace*{-0.02in}
\end{bmatrix}  \hspace*{-0.02in} - \hspace*{-0.02in} \eta I \geq 0,
\label{eq:polya_blockdiag_simplex}
\end{align}
\end{small}
\hspace*{-0.15in} where we have denoted the elements of $I_{d_p+e}$ by $g^{(i)} \in \mathbb{N}^l, i=1,\cdots,L$ and have denoted the elements of $I_{d_p+d_a+e}$ by $q^{(i)}, i=1,\cdots,M$. For any  $\gamma \in \mathbb{R}$, if there exist $e \geq 0$ such that SDP~\eqref{eq:polya_blockdiag_simplex} is feasible, then $\gamma \leq \gamma^*$. If for a positive $\gamma$, there exist $e \geq 0$ such that SDP~\eqref{eq:polya_blockdiag_simplex} has a solution $P_h, h \in I_{d_p}$, then $V=x^T \left(  \sum_{h \in I_{d_p}} P_h \alpha_1^{h_1} \cdots \alpha_l^{h_l} \right)x$ is a Lyapunov function proving stability of $\dot{x}(t)=A(\alpha)x(t), \alpha \in \Delta^l$.  See~\cite{reza_peet_tac2013} for a complexity analysis on SDP~\eqref{eq:polya_blockdiag_simplex}.


SDPs such as~\eqref{eq:polya_blockdiag_simplex} can be solved in polynomial time using interior-point algorithms such as the central path primal-dual algorithms in~\cite{monteiro,helmberg,alizadeh}.
Fortunately, Problem~\eqref{eq:polya_blockdiag_simplex} has block-diagonal structure.
Block-diagonal structure in SDP constraints can be used to design massively parallel algorithms, an approach which was applied to Problem~\eqref{eq:polya_blockdiag_simplex} in~\cite{reza_peet_tac2013}.

\noindent \textbf{\\Numerical Example 2:}
Consider the system $\dot{x}(t) = A(\alpha) x(t)$, where 
\[
A(\alpha)= A_1 \alpha_1^3 + A_2 \alpha_1^2 \alpha_2 + A_3 \alpha_1 \alpha_3^2 + A_4 \alpha_1 \alpha_2 \alpha_3 + A_5 \alpha_2^3+ A_6 \alpha_3^3,
\]
where 
\[
\alpha \in T_L := \{\alpha \in \mathbb{R}^3: \sum_{i=1}^3 \alpha_i=2L+1, L \leq \alpha_i \leq 1 \}
\]
 and
\begin{small}
\begin{align*}
&A_1=\begin{bmatrix}
-0.57 & -0.44 & 0.33 & -0.07 \\
-0.48 & -0.60 & 0.30 &  0 \\
-0.22 & -1.12 & 0.08 & -0.24 \\
 1.51 & -0.42 & 0.67 & -1.00
\end{bmatrix} 
\quad A_2=\begin{bmatrix}
-0.09 & -0.16  &  0.3  &  -1.13\\
-0.15 &  -0.17 & -0.02 &   0.82\\
0.14  &   0.06 &  0.02 & -1\\
0.488 &  0.32  &  0.97 & -0.71
\end{bmatrix} \\
&A_3=\begin{bmatrix}
-0.70 & -0.29 & -0.18 &  0.31\\
0.41  & -0.76 & -0.30 & -0.12\\
-0.05 &  0.35 & -0.59 &  0.91\\
1.64  &  0.82 &  0.01 &    -1
\end{bmatrix}
\; A_4=\begin{bmatrix}
0.72 &  0.34 & -0.64 &  0.31\\
-0.21 & -0.51 &  0.59 &  0.07\\
 0.27 &  0.49 & -0.84 & -0.94\\
-1.89 & -0.66 &  0.27 &  0.41
\end{bmatrix}   \\
&A_5=\begin{bmatrix}
-0.51 & -0.47 & -1.38 &  0.17\\
1.18  & -0.62 & -0.29 &  0.35\\
-0.65 &  0.01 & -1.44 & -0.04\\
-0.74 & -1.22 &  0.60 & -1.47\\
\end{bmatrix}      
\;A_6=\begin{bmatrix}
-0.201 & -0.19  &  -0.55 &  0.07\\
0.803  & -0.42  & -0.20  &  0.24\\
-0.440 &  0.01  & -0.98 & -0.03\\
0      & -0.83  &  0.41 & -1
\end{bmatrix} .    
\end{align*}
\end{small}
We would like to solve the following optimization problem.
\begin{align}
& L^*= \min \;\;\; L \nonumber \\
& \text{subject to } \;\; \dot{x}(t) = A(\alpha) x(t) \text{ is stable for all } \alpha \in   T_L.  
\label{eq:optim1}
\end{align}

\begin{figure}[htbp]
\centering
\includegraphics[scale=0.27]{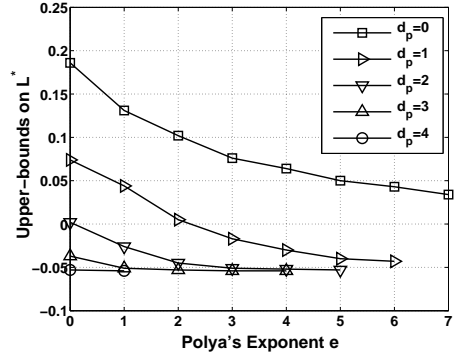}
\caption{Upper-bounds for $L^*$ in Problem~\eqref{eq:optim1} for different Polya's exponents $e$ and different degrees of $P$}
\label{fig:simplex}
\end{figure}

We first represent $T_L$ using the unit simplex $\Delta^3$ as follows. Define the map $f: \Delta^3 \rightarrow T_L$ as
\[
f(\alpha) = [f_1(\alpha), f_2(\alpha), f_3(\alpha)],
\]
where $f_i(\alpha) = 2 \vert L \vert (\alpha_i - 0.5)$.
Then, we have $\{A(\alpha): \alpha \in T_L \} = \{ A(f(\beta)), \beta \in \Delta^3 \}$. Thus, the following optimization problem is equivalent to Problem~\eqref{eq:optim1}. 
\begin{align}
& L^*= \min \;\;\; L \nonumber \\
& \text{subject to } \;\; \dot{x}(t) = A(f(\beta)) x(t) \text{ is stable for all } \beta \in   \Delta^3.  
\label{eq:optim2}
\end{align}
We solved Problem~\eqref{eq:optim2} using bisection on $L$. For each $L$, we used Theorem~\ref{thm:polya_simplex} to verify the inequality in~\eqref{eq:OOP_simplex} using Polya's exponents $e=1$ to 7 and $d_p=1$ to 4 as degrees of $P$. Figure~\ref{fig:simplex} shows the computed upper-bounds on $L^*$ for different $e$ and $d_p$. The best upper-bound found by the algorithm is $-0.0504$.

For comparison, we solved the same problem using SOSTOOLS~\cite{sostools2013} and Putinar's Positivstellensatz (see Case 4 of Section~\ref{sec:algorithms}). By computing a Lyapunov function of degree two in $x$ and degree one in $\beta$, SOSTOOLS certified $L=-0.0504$ as an upper-bound for $L^*$. This is the same as the upper-bound computed by Polya's algorithm. The CPU time required for SOSTOOLS to compute the upper-bound on a Core i7 machine with 64 GB of RAM was 22.3 minutes, whereas the Polya's algorithm only required 7.1 seconds to compute the same upper-bound.

\noindent \textbf{\\Case 3. $A(\alpha)$ is a polynomial with $\alpha \in \Phi^l$:}

Given $A_h \in \mathbb{R}^{n \times n}$ for $h \in E_d$ as defined in~\eqref{eq:Ed}, we address the problem of stability analysis of System~\eqref{eq:sys_lin} with $A(\alpha)=\sum_{h \in E_d} A_h \alpha_1^{h_1} \cdots \alpha_l^{h_l} $ for all $\alpha \in \Phi^l:= \{ x \in \mathbb{R}^{n}: \vert x_i \vert \leq r_i \}$.
 From Lyapunov theory, System~\eqref{eq:sys_lin} is asymptotically stable for all $\alpha \in \Phi^l$ if and only if $\gamma^* > 0$ in the following optimization of polynomials problem.
\begin{align}
&\gamma^*=\max_{\gamma \in \mathbb{R}, P \in \mathbb{R}[\alpha]} \;\; \gamma \nonumber \\
& \text{subject to } \begin{bmatrix}
P(\alpha)      &  0           \\ 
0      & -A(\alpha)^T P(\alpha) - P(\alpha) A(\alpha) 
\label{eq:OOP_multisimplex} 
\end{bmatrix}
-
\gamma I \geq 0  \text{ for all } \alpha \in \Phi^l
\end{align}
As in Case 2 of Section~\ref{sec:algorithms}, by applying bisection algorithm on $\gamma$ and using a multi-simplex version of Polya's algorithm (such as Theorem~\ref{thm:polya_multi-simplex}) as a test for feasibility of Constraint~\eqref{eq:OOP_multisimplex} we can compute lower bounds on $\gamma^*$.
Suppose there exists a matrix-valued multi-homogeneous polynomial (defined in~\eqref{eq:multi-homog_poly}) $Q$ of degree vector $d_q \in \mathbb{N}^l$ ($d_{q_i}$ is the degree of $\beta_i$) such that
\begin{equation}
\{P(\alpha) \in \mathbb{S}^n: \alpha \in \Phi^l \} = \{ Q(\beta,\eta) \in \mathbb{S}^n: \beta, \eta \in \mathbb{R}^l \text{ and } (\beta_i,\eta_i) \in \Delta^2 \text{ for } i=1,\cdots,l \}.
\label{eq:PQ_sets}
\end{equation}
Likewise, suppose there exists a matrix-valued multi-homogeneous polynomial $B$ of degree vector $d_b \in \mathbb{N}^l$ ($d_{b_i}$ is the degree of $\beta_i$) such that
\[
\{ A(\alpha) \in \mathbb{S}^n: \alpha \in \Phi^l \} = \{ B(\beta,\eta) \in \mathbb{S}^n: \beta,\eta \in \mathbb{R}^l  \text{ and } (\beta_i,\eta_i) \in \Delta^2 \text{ for } i=1, \cdots,l \}.
\]
Given $\gamma \in \mathbb{R}$, it follows from Theorem~\ref{thm:polya_multi-simplex} that the inequality condition in~\eqref{eq:OOP_multisimplex} holds for all $\alpha \in \Phi^l$ if there exist $e \geq 0$ such that
\begin{align}
& \hspace*{-0.2in} \left( \prod_{i=1}^l (\beta_i+\eta_i)^e \right) \left( Q(\beta,\eta) - \gamma I \left( \prod_{i=1}^l (\beta_i + \eta_i)^{d_{p_i}} \right) \right)  \label{eq:polya_ineq3} \\
\text{ and } \nonumber \\
  &  \hspace*{-0.35in} -\left( \prod_{i=1}^l (\beta_i + \eta_i)^e \right) \left(B^T(\alpha,\beta)Q(\beta,\eta)+Q(\beta,\eta)B(\beta,\eta) + \gamma I \left( \prod_{i=1}^l (\beta_i + \eta_i)^{d_{pa_i}} \right) \right),
 \label{eq:polya_ineq4}
\end{align}
have all positive coefficients where $d_{p_i}$ is the degree of $\alpha_i$ in $P(\alpha)$ and $d_{pa_i}$ is the degree of $\alpha_i$ in $P(\alpha)A(\alpha)$. Suppose $Q$ and $B$ are of the forms
\begin{equation}
Q(\beta,\eta) = \sum_{ \substack{ h,g \in \mathbb{N}^l \\ h+g=d_q }} Q_{h,g} \beta_1^{h_1} \eta_1^{g_1} \cdots \beta_l^{h_l} \eta_l^{g_l}
\label{eq:Q_multisimplex}
\end{equation}
and 
\begin{equation}
B(\beta,\eta) = \sum_{ \substack{ h,g \in \mathbb{N}^l \\ h+g=d_b }} B_{h,g} \beta_1^{h_1} \eta_1^{g_1} \cdots \beta_l^{h_l} \eta_l^{g_l}.
\label{eq:B_multisimplex}
\end{equation}
By combining~\eqref{eq:Q_multisimplex} and~\eqref{eq:B_multisimplex} with~\eqref{eq:polya_ineq3} and~\eqref{eq:polya_ineq4} we find that for a given $\gamma \in \mathbb{R}$, the inequality condition in~\eqref{eq:OOP_multisimplex} holds for all $\alpha \in \Phi^l$ if there exist some $e \geq 0$ such that
\begin{align}
& \sum_{\substack{ h,g \in \mathbb{N}^l \\ h+g=d_q }} \hspace*{-0.06in}  f_{\{q,r\}, \{h,g \}} Q_{h,g} > 0 \;\; \text{ for all } \; q,r \in \mathbb{N}^l: q+r=d_q+e \cdot \mathbf{1} \;\; \text{ and} \nonumber \\
& \sum_{\substack{ h,g \in \mathbb{N}^l \\ h+g=d_q }} \hspace*{-0.08in} M_{\{s,t\},\{h,g\}}^T Q_{h,g} + Q_{h,g} M_{\{s,t\},\{h,g\}} \hspace*{-0.03in} < 0  \text{ for all } s,t \in \mathbb{N}^l \hspace*{-0.03in} : \hspace*{-0.03in} s+t=d_q+d_b+e \cdot \mathbf{1},
\label{eq:polya_multisimplex_conditions}
\end{align}
where $\mathbf{1} \in \mathbb{N}^l $ is the vector of ones and where we define $f_{\{q,r\}, \{h,g \}} \in \mathbb{R}$ to be the coefficient of $Q_{h,g} \beta^q \eta^r$ after expansion of~\eqref{eq:polya_ineq3}. Likewise, we define $M_{\{s,t\},\{h,g\}} \in \mathbb{R}^{n \times n}$ to be the coefficient of $Q_{h,g}\beta^s \eta^t$ after expansion of~\eqref{eq:polya_ineq4}.
See~\cite{reza_CDC_hypercube} for recursive formulae for calculating $f_{\{q,r\}, \{h,g \}}$ and $M_{\{s,t\},\{h,g\}}$. Similar to Case 2, Conditions~\eqref{eq:polya_multisimplex_conditions} are an SDP (See~\cite{reza_CDC_hypercube} for a complexity analysis on this SDP). For any $\gamma \in \mathbb{R}$, if there exist $e \geq 0$ and $\{ Q_{h,g} \}$ that satisfy~\eqref{eq:polya_multisimplex_conditions}, then $\gamma \leq \gamma^*$ as defined in~\eqref{eq:OOP_multisimplex}. Furthermore, if  $\gamma$ is positive, then $\dot{x}(t)=A(\alpha)x(t)$ is asymptotically stable for all $\alpha \in \Phi^l$.

\noindent\textbf{\\Numerical Example 3a:} 
Consider the system $\dot{x}(t) = A(\alpha) x(t)$, where 
\[
A(\alpha)= A_0 + A_1 \alpha_1^2 + A_2 \alpha_1 \alpha_2 \alpha_3 + A_3 \alpha_1^2 \alpha_2 \alpha_3^2,  
\]
\[
\alpha_1 \in [-1,1], \, \alpha_2 \in [-0.5,0.5], \, \alpha_3 \in [-0.1,0.1],  
\] 
where 
\begin{small}
\begin{align*}
&A_0=\begin{bmatrix}
   -3.0  &     0 &  -1.7 &   3.0\\
   -0.2  & -2.9 & -1.7  & -2.6\\
    0.6  &  2.6 &  -5.8  & -2.6\\
   -0.7  &  2.9  & -3.3  & -2.1
\end{bmatrix} 
A_1=\begin{bmatrix}
    2.2 &  -5.4 &  -0.8  & -2.2\\
    4.4 &  1.4  & -3.0   & 0.8\\
   -2.4 &  -2.2 &   1.4  &  6.0\\
   -2.4 &  -4.4 &  -6.4  &  0.18
\end{bmatrix} \\
&A_2=\begin{bmatrix}
   -8.0 & -13.5 &  -0.5 &  -3.0\\
   18.0 &  -2.0  &  0.5 & -11.5\\
    5.5 & -10.0 &   3.5  &  9.0\\
   13.0 &   7.5 &   5.0 &  -4.0
\end{bmatrix}
A_3=\begin{bmatrix}
    3.0  &  7.5 &   2.5 &  -8.0\\
    1.0  &  0.5 &   1.0 &   1.5\\
   -0.5  & -1.0 &   1.0 &   6.0\\
   -2.5  & -6.0 &   8.5 &  14.25
\end{bmatrix}.    
\end{align*}
\end{small}
\hspace*{-0.1in} The problem is to investigate asymptotic stability of this system for all $\alpha$ in the given intervals using the method in Case 3 of Section~\ref{sec:applications_lin}.
We first represented $A(\alpha)$ over $[-1,1] \times [-0.5,0.5] \times [-0.1,0.1]$ by a multi-homogeneous polynomial $B(\beta,\eta)$ with $(\beta_i,\eta_i) \in \Delta^2$ and with the degree vector $d_b=[2,1,2]$ (see~\cite{reza_CDC_hypercube} for an algorithm which finds $B$ and see Case 2 of Section~\eqref{sec:algorithms} for an example). Then, by applying Theorem~\ref{thm:polya_multi-simplex} (as in~\eqref{eq:polya_ineq3} and~\eqref{eq:polya_ineq4}) with $\gamma =0.1,e=1$ and $d_p=[1,1,1]$, we set-up the inequalities in~\eqref{eq:polya_multisimplex_conditions} with $d_q=[1,1,1]$. By using semi-definite programming, we solved the inequalities and computed the following Lyapunov function as a certificate for asymptotic stability of $\dot{x}(t) = A(\alpha) x(t)$ for all $\alpha_1 \in [-1,1], \, \alpha_2 \in [-0.5,0.5], \, \alpha_3 \in [-0.1,0.1]$.
\begin{align*}
V(x,\beta,\eta)= x^T Q(\beta,\eta)x  = x^T ( & \beta_1 (Q_1 \beta_2 \beta_3 + Q_2 \beta_2 \eta_3 +Q_3 \eta_2 \beta_3 + Q_4  \eta_2 \eta_3)  \\
 + & \eta_1 (Q_5 \beta_2 \beta_3 + Q_6 \beta_2 \eta_3 + Q_7 \eta_2 \beta_3 + Q_8  \eta_2 \eta_3 ) )x,
\end{align*}
where $\beta_1= 0.5 \alpha_1 +0.5, \beta_2=\alpha_2+0.5, \beta_3=5\alpha_3+0.5, \eta_1= 1- \beta_1, \eta_2 = 1- \beta_2, \eta_3 = 1- \beta_3$ and
\begin{small}
\begin{align*}
&Q_1=\begin{bmatrix}
   5.807  &     0.010 &  -0.187 &  -1.186\\
    0.010  &  5.042 & -0.369  &  0.227\\
    -0.187  & -0.369 &  8.227  & -1.824\\
   -1.186  &   0.227  & -1.824  &  8.127
\end{bmatrix} 
Q_2=\begin{bmatrix}
    7.409 &  -0.803 &  1.804  & -1.594\\
    -0.803 &  6.016  &  0.042   & -0.538\\
   1.804 &   0.042 &   7.894  &  -1.118\\
   -1.594 &  -0.538 &  -1.118  &  8.590
\end{bmatrix} \\
&Q_3=\begin{bmatrix}
    6.095 & -0.873 &   0.512 &  -1.125\\
   -0.873 &  5.934  &   -0.161 &  0.503\\
    0.512 &  -0.161 &   7.417  &  -0.538\\
   -1.125 &    0.503 &   -0.538 &  6.896
\end{bmatrix}
Q_4=\begin{bmatrix}
    5.388  &  0.130 &   -0.363 &  -0.333\\
    0.130  &  5.044 &   -0.113 &   -0.117\\
   -0.363  & -0.113 &    6.156 &   -0.236\\
   -0.333  & -0.117 &   -0.236 &   5.653
\end{bmatrix}   \\
&Q_5=\begin{bmatrix}
   7.410  &  -0.803 &  1.804 &  -1.594\\
   -0.803  &  6.016 &  0.042  &  -0.538\\
    1.804  &  0.042 &  7.894  & -1.118\\
   -1.594  &  -0.538  & -1.118  &  8.590
\end{bmatrix} 
Q_6=\begin{bmatrix}
    5.807 &   0.010 &  -0.187  & -1.186\\
     0.010 &   5.042  & -0.369   &  0.227\\
   -0.187 &  -0.369 &    8.227  & -1.824\\
   -1.186 &   0.227 &  -1.824  &   8.127
\end{bmatrix} \\
&Q_7=\begin{bmatrix}
    5.388 & 0.130 &   -0.363 &  -0.333\\
   0.130 &  5.044  & -0.113 & -0.117\\
   -0.363 & -0.113 &  6.156  &   -0.236\\
   -0.333 &  -0.117 &   -0.236 &  5.653
\end{bmatrix}
Q_8=\begin{bmatrix}
    6.095  &   -0.873 &   0.512 &  -1.125\\
     -0.873  &  5.934 &  -0.161 &   0.503\\
   0.512  & -0.161 &   7.417 &   -0.538\\
    -1.125  &  0.503 &  -0.538 &   6.896
\end{bmatrix}    .
\end{align*}
\end{small}

\noindent\textbf{\\Numerical Example 3b:} In this example, we used the same method as in Example 3a to find lower bounds on $r^*= \max \, r$ such that $\dot{x}(t) = A(\alpha) x(t)$ with
\[
A(\alpha) = A_0+\sum_{i=1}^4 A_i \alpha_i,
\]
\begin{small}
\begin{align*}
&A_0=\begin{bmatrix}
   -3.0  &  0   & -1.7 &   3.0\\
   -0.2  & -2.9 & -1.7  & -2.6\\
    0.6  &  2.6 & -5.8  & -2.6\\
   -0.7  &  2.9 & -3.3  & -2.4
\end{bmatrix} 
A_1=\begin{bmatrix}
    1.1 &  -2.7 &  -0.4  & -1.1\\
    2.2 &  0.7  &  -1.5   & 0.4\\
   -1.2 &  -1.1 &   0.7  &  3.0\\
   -1.2 &  -2.2 &  -3.2  &  -1.4
\end{bmatrix} \\
&A_2=\begin{bmatrix}
   1.6 & 2.7 &  0.1 &  0.6\\
   -3.6 &  0.4  &  -0.1 & 2.3\\
    -1.1 & 2 &  -0.7  &  -1.8\\
   -2.6 &  -1.5 &  -1.0 &  0.8
\end{bmatrix}
A_3=\begin{bmatrix}
    -0.6  & 1.5 &  0.5 & -1.6 \\
     0.2 & -0.1 & 0.2  &  0.3 \\
     -0.1 & -0.2  & -0.2  & 1.2  \\
    -0.5 & -1.2 & 1.7  & -0.1 
\end{bmatrix} \\
& \qquad \qquad \qquad \qquad A_4=\begin{bmatrix}
    -0.4  & -0.1 &  -0.3 & 0.1 \\
     0.1 & 0.3 & 0.2  &  0.0 \\
     0.0 & 0.2  & -0.3  & 0.1  \\
     0.1 & -0.2 & -0.2  & 0.0 
\end{bmatrix}.     
\end{align*}
\end{small}
\hspace*{-0.09in} is asymptotically stable for all $\alpha \in \{ \alpha \in \mathbb{R}^4: \vert \alpha_i \vert \leq r \}$. Table~\ref{tab:example3b} shows the computed lower bounds on $r^*$ for different degree vectors $d_q$ (degree vector of $Q$ in~\eqref{eq:PQ_sets}). In all of the cases, we set the Polya's exponent $e=0$.
For comparison, we have also included the lower-bounds computed by the methods of~\cite{bliman_SIAM_2004} and~\cite{chesi_hypercube_2005} in Table~\ref{tab:example3b}. 

\begin{scriptsize}
\renewcommand{\tabcolsep}{2.5pt} 
\begin{table}[h]
\caption{The lower-bounds on $r^*$ computed by the method in Case 3 of Section~\ref{sec:applications_lin} and methods in~\cite{bliman_SIAM_2004} and~\cite{chesi_hypercube_2005} -  $i^{th}$ entry of $d_q$ is the degree of $\beta_i$ in~\eqref{eq:PQ_sets}}
\label{tab:example3b}
\begin{center}
\begin{tabular}{ c|c|c|c|c|c|c|c|c}
\cline{2-8}
& $d_q=$[0,0,0,0] & $d_q=$[0,1,0,1] & $d_q=$[1,0,1,0] & $d_q=$[1,1,1,1] &  $d_q=$[2,2,2,2] & Ref.\cite{bliman_SIAM_2004} & Ref.\cite{chesi_hypercube_2005} \\ 
\hline
\hspace*{-5.5pt} \vline bound on $r^*$ & 0.494 & 0.508 & 0.615 & 0.731 & 0.840 & 0.4494 & 0.8739 \\ 
\hline 
\end{tabular} 
\end{center}
\end{table}
\end{scriptsize}


\subsection{Nonlinear stability analysis}
Consider nonlinear systems of the form
\begin{equation}
\dot{x}(t)=f(x(t)),
\label{eq:sys_nonlin}
\end{equation}
where $f:\mathbb{R}^n \rightarrow \mathbb{R}^n$ is a degree $d_f$ polynomial. Suppose the origin is an isolated equilibrium of~\eqref{eq:sys_nonlin}. In this section, we address local stability of the origin in the following sense.

\begin{lemma}
Consider the System~\eqref{eq:sys_nonlin} and let $Q \subset \mathbb{R}^{n \times n}$ be an open set containing the origin. Suppose there exists a continuously differentiable function $V$ which satisfies
\begin{equation}
V(x) > 0 \text{ for all } x \in Q \setminus \{0\}, V(0)=0
\label{eq:V}
\end{equation}
and 
\begin{equation}
\langle \nabla V,f(x) \rangle < 0  \text{ for all } x \in Q \setminus \{0\}.
\label{eq:Vdot}
\end{equation}
Then the origin is an asymptotically stable equilibrium of System~\eqref{eq:sys_nonlin}, meaning that for every
$x(0) \in \{ x \in \mathbb{R}^n : \{ y:V(y) \leq V(x) \} \subset Q\}$,  $ \lim_{t \to \infty} x(t) = 0$.
\label{lemma:Lyap}
\end{lemma}

 Since existence of polynomial Lyapunov functions is necessary and sufficient for stability of~\eqref{eq:sys_nonlin} on any compact set~\cite{peet_papa_TAC2012}, we can formulate the problem of stability analysis of~\eqref{eq:sys_nonlin} as follows.
\begin{align}
& \gamma^* = \max_{\gamma, c_\beta \in \mathbb{R}} \;\; \gamma \nonumber \\
& \text{subject to} \begin{bmatrix}
\sum_{\beta \in E_d} c_\beta x^\beta - \gamma x^Tx &  0\\ 
0 & \hspace*{-0.1in}  -\langle \nabla \sum_{\beta \in E_d} c_\beta x^\beta, f(x) \rangle - \gamma x^Tx 
\end{bmatrix}  \geq 0 \text{ for all } x \in Q.
\label{eq:OOP_nonlin0}
\end{align}
Conditions~\eqref{eq:V} and~\eqref{eq:Vdot} hold if and only if there exist $d \in \mathbb{N}$ such that $\gamma^* > 0$.
In Sections~\ref{sec:handelman_nonlin} and~\ref{sec:polya_nonlin}, we discuss two alternatives to SOS programming for solving~\eqref{eq:OOP_nonlin0}. These methods apply Polya's theorem and Handelman's theorem to Problem~\eqref{eq:OOP_nonlin0} (as described in Cases 2 and 3 in Section~\ref{sec:algorithms}) to find lower bounds on $\gamma^*$. 
See~\cite{Handelman_Sankaranarayanan} for a different application of Handelman's theorem and intervals method in nonlinear stability. 
Also, see~\cite{CPA_sigurdur} for a method of computing continuous piecewise affine Lyapunov functions using linear programming and a \textit{triangulation} scheme for polytopes.

\subsubsection{Application of Handelman's theorem in nonlinear stability analysis}
\label{sec:handelman_nonlin}
Recall that every convex polytope can be represented as
\begin{equation}
\Gamma^K:= \{ x \in \mathbb{R}^n: w_i^T x + u_i \geq 0, i=1, \cdots, K \}
\label{eq:polytope}
\end{equation}
for some $w_i \in \mathbb{R}^n$ and $u_i \in \mathbb{R}$.
Suppose $\Gamma^K$ is bounded and the origin is in the interior of $\Gamma^K$. In this section, we would like to investigate asymptotic stability of the equilibrium of System~\eqref{eq:sys_nonlin} by verifying positivity of $\gamma^*$ in Problem~\eqref{eq:OOP_nonlin0} with $Q=\Gamma^K$.

Unfortunately, Handelman's theorem (Theorem~\ref{thm:Handelman}) does not parameterize polynomials which have zeros in the interior of a given polytope. To see this, suppose a polynomial $g$ ($g$ is not identically zero) is zero at $x=a$, where $a$ is in the interior of a polytope $\Gamma^K:= \{ x \in \mathbb{R}^n: w_i^T x + u_i \geq 0, i=1, \cdots, K \}$. Suppose there exist $b_\alpha \geq 0, \alpha \in \mathbb{N}^K$ such that for some $d \in \mathbb{N}$,
\[
g(x) = \sum _{\substack{\alpha \in \mathbb{N}^K \\ \Vert \alpha_i \Vert_1 \leq d}} b_\alpha (w_1^T x+u_1)^{\alpha_1} \cdots (w_K^T x+u_K)^{\alpha_K}.
\]
Then, 
\[
g(a) = \sum _{\substack{\alpha \in \mathbb{N}^K \\ \Vert \alpha_i \Vert_1 \leq d}} b_\alpha (w_1^T a+u_1)^{\alpha_1} \cdots (w_K^T a+u_K)^{\alpha_K} = 0.
\]
 From the assumption $a \in \text{int}(\Gamma^K)$ it follows that $w_i^Ta+u_i > 0$ for $i=1, \cdots,K$. Hence $b_\alpha < 0$ for at least one $\alpha \in \{ \alpha \in \mathbb{N}^K : \Vert \alpha \Vert_1 \leq d \} $. This contradicts with the assumption that all $b_\alpha \geq 0$.

Based on the above reasoning, one cannot readily use Handelman's theorem to verify the Lyapunov inequalities in~\eqref{eq:V}. In~\cite{reza_CDC_2014}, a combination of Handelman's theorem and a decomposition scheme was applied to Problem~\eqref{eq:OOP_nonlin0} with $Q=\Gamma^K$. Here we outline this result. First, consider the following definitions.

\begin{definition}
Given a bounded polytope of the form $\Gamma^K:= \{ x \in \mathbb{R}^n: w_i^T x + u_i \geq 0, i=1, \cdots, K \}$, we call 
\[
\zeta^i(\Gamma^K):=\{ x \in \mathbb{R}^n: w_i^Tx+u_i = 0 \text{ and } w_j^Tx+u_j \geq 0 \text{ for } j \in \{1, \cdots,K \} \}
\]
the $i-$th facet of $\Gamma^K$ if $\zeta^i(\Gamma^K) \neq \emptyset$.
\end{definition}

\begin{definition}
Given a bounded polytope of the form $\Gamma^K:= \{ x \in \mathbb{R}^n: w_i^T x + u_i \geq 0, i=1, \cdots, K \}$, we call $D_\Gamma:=\{D_i\}_{i=1,\cdots,L}$ a $D-$\textit{decomposition} of $\Gamma^K$ if 
\begin{equation*}
D_i := \{ x \in \mathbb{R}^n: h^T_{i,j} x + g_{i,j} \geq 0, j=1, \cdots,m_i \} \;\; \text{ for some } h_{i,j} \in \mathbb{R}^n, \, g_{i,j} \in \mathbb{R}
\label{eq:decomposition}
\end{equation*}
such that $\cup_{i=1}^L D_i = \Gamma^K,\, \cap_{i=1}^L D_i = \{0\}$ and $\text{int}(D_i) \cap \text{int}(D_j) = \emptyset$. 
\end{definition}

Consider System~\eqref{eq:sys_nonlin} with polynomial $f$ of degree $d_f$. Given $w_i, h_{i,j} \in \mathbb{R}^n$ and $u_i, g_{i,j} \in \mathbb{R}$, let $\Gamma^K:= \{ x \in \mathbb{R}^n: w_i^T x + u_i \geq 0, i=1, \cdots, K \}$ with $D-$decomposition $D_\Gamma:=\{D_i\}_{i=1,\cdots,L}$, where $D_i := \{ x \in \mathbb{R}^n: h^T_{i,j} x +g_{i,j} \geq 0, j=1, \cdots,m_i \}$.
Let us denote the elements of the set
\[
E_{d,n}:=\{ \alpha \in \mathbb{N}^n : \sum_{i=1}^n \alpha_i \leq d \}
\]
by $\lambda^{(k)}, k=1, \cdots, B$, where $B$ is the cardinality of $E_{d,e}$. For any $\lambda^{(k)}$, let $p_{\{ \lambda^{(k)}, \alpha,i \}}$ be the coefficient of $b_{i, \alpha} \, x^{\lambda^{(k)}}$ in
\begin{equation}
P_i(x) := \sum_{\alpha \in E_{d,m_i}} b_{i,\alpha} \prod_{j=1}^{m_i} (h_{i,j}^T x + g_{i,j})^{\alpha_j}, \; x \in \mathbb{R}^n, b_{i,\alpha} \in \mathbb{R}.
\label{eq:P_i}
\end{equation}
Let us denote the cardinality of $E_{d,m_i}$ by $N_i$ and denote the vector of all coefficients $b_{i, \alpha}$ by $b_i \in \mathbb{R}^{N_i}$. Define the map  $F_{i}: \mathbb{R}^{N_i} \times \mathbb{N} \rightarrow \mathbb{R}^{B}$ as 
\[
F_{i}(b_i,d) := \left[ 
 \sum_{\substack{ \alpha \in E_{d,m_i} }}  p_{\{\lambda^{(1)},\alpha,i\}} b_{i,\alpha} \; , \; \cdots \; , \;  \sum_{\substack{ \alpha \in E_{d,m_i} }}  p_{\{\lambda^{(
B)},\alpha,i\}} b_{i,\alpha} \right]^T 
\]
for  $i=1, \cdots,L$. In other words, $F_i(b_i,d)$ is the vector of coefficients of $P_i(x)$ after expansion.

\noindent Define $H_{i}: \mathbb{R}^{N_i} \times \mathbb{N} \rightarrow \mathbb{R}^{Q} $ as
\[
H_{i}(b_i,d) := \left[
 \sum_{\substack{\alpha \in E_{d,m_i} }}  p_{\{\delta^{(1)},\alpha,i\}} b_{i,\alpha} \; , \; \cdots \; , \;  \sum_{\substack{ \alpha \in E_{d,m_i} }}  p_{\{\delta^{(
Q)},\alpha,i\}} b_{i,\alpha} \right]^T
\]
for $i = 1, \cdots,L$, where we have denoted the elements of $\{ \delta \in \mathbb{N}^n: \delta = 2e_j \text{ for } j=1, \cdots,n \}$ by $\delta^{(k)}, k=1, \cdots, Q$, where $e_j$ are the canonical basis for $\mathbb{N}^n$. In other words, $H_{i}(b_i,d)$ is the vector of coefficients of square terms of $P_i(x)$ after expansion.     \\

\noindent Define $R_i(b_i,d): \mathbb{R}^{N_i} \times \mathbb{N} \rightarrow \mathbb{R}^C$ as
\[
R_i(b_i,d):=
\left[ b_{i,\beta^{(1)}} \; , \; \cdots \; , \; b_{i,\beta^{(C)}} \right]^T,
\]
for $i = 1, \cdots,L$, where we have denoted the elements of the set 
\[
\{ \beta \in E_{d,m_i} :  \beta_j=0 \text{ for } j \in \{ j \in \mathbb{N}: g_{i,j} = 0 \} \}   
\]
 by $\beta^{(k)}$ for $k=1,\cdots,C$. \\

\noindent Define $J_{i}: \mathbb{R}^{N_i} \times \mathbb{N} \times \{1, \cdots, m_i \} \rightarrow \mathbb{R}^{B} $ as
\[
J_{i}(b_i,d,k) := \left[ 
 \sum_{\substack{ \alpha \in E_{d,m_i} \\ \alpha_k=0} }  p_{\{\lambda^{(1)},\alpha,i\}} b_{i,\alpha} \; , \; \cdots \; , \;  \sum_{\substack{ \alpha \in E_{d,m_i} \\ \alpha_k=0 }}  p_{\{\lambda^{(B)},\alpha,i\}} b_{i,\alpha} \right]^T
\]
for $i,k = 1, \cdots,L$. In other words, $J_{i}(b_i,d,k)$ is the vector of coefficients of restriction of $P_i$ to $h_{i,k}^Tx + g_{i,k} = 0$, after expansion. \\

\noindent  Define $G_{i}: \mathbb{R}^{N_i} \times \mathbb{N} \rightarrow \mathbb{R}^{Z} $ as
\[
G_{i}(b_i,d) := \left[ 
 \sum_{\substack{ \alpha \in E_{d,m_i}  }}  q_{\{\eta^{(1)},\alpha,i\}} b_{i,\alpha} \; , \; \cdots \; , \;  \sum_{\substack{ \alpha \in E_{d,m_i}  }}  q_{\{\eta^{(P)},\alpha,i\}} b_{i,\alpha} \right]^T
\]
for $i=1, \cdots,L$, where we have denoted the elements of $E_{d+d_f-1,n}$ by $\eta^{(k)}$ for $k=1, \cdots, Z$. For any $\eta^{(k)} \in E_{d+d_f-1,n}$, we define $q_{\{\eta^{(k)},\alpha,i\}}$ as the coefficient of $b_{i,\alpha} \, x^{\eta^{(k)}}$ in $\langle \nabla P_i(x),f(x) \rangle$, where $P_i(x)$ is defined in~\eqref{eq:P_i}.

\noindent Finally, given $i,j \in \{1, \cdots,L \}, i \neq j$, let
\begin{small}
\[
\Lambda_{i,j} := \left\lbrace   k,l \in \mathbb{N} : k \in \{ 1, \cdots, m_i \}, l \in \{ 1, \cdots, m_j \}: \zeta^k(D_i) \neq \emptyset \text{ and }  \zeta^k(D_i) = \zeta^l(D_j) \right\rbrace 
\]
\end{small}
If there exist $d \in \mathbb{N}$ such that $\max \gamma$ in the linear program 
\begin{align}
&\max_{\gamma \in \mathbb{R} , b_i \in \mathbb{R}^{N_i}, c_i \in \mathbb{R}^{M_i}} \;\;\; \gamma  \nonumber  \\
&  \text{subject to } \;\; b_{i} \geq \mathbf{0}   && \text{ for } i=1, \cdots, L \nonumber \\
& \hspace*{0.73in} c_{i} \leq \mathbf{0}           && \text{ for } i=1, \cdots, L \nonumber  \\
& \hspace*{0.73in} R_i(b_{i},d) = \mathbf{0}       && \text{ for } i=1, \cdots, L        \nonumber \\
& \hspace*{0.73in} H_i(b_i,d) \geq \gamma \cdot \mathbf{1}      &&  \text{ for } i=1, \cdots, L \nonumber \\
& \hspace*{0.73in} H_i(c_i,d+d_f-1) \leq -\gamma \cdot \mathbf{1}  && \text{ for } i=1, \cdots, L \nonumber \\
& \hspace*{0.73in} G_i(b_i,d) = F_i(c_i,d+d_f-1)   &&\text{ for } i=1, \cdots, L \nonumber \\
& \hspace*{0.73in} J_{i}(b_i,d,k) = J_{j}(b_j,d,l)  &&  \text{ for } i,j=1, \cdots, L \text{ and } k,l \in \Lambda_{i,j}
\label{eq:lin_prog}
\end{align}
is positive, then the origin is an asymptotically stable equilibrium for System~\eqref{eq:sys_nonlin} and
\[
V(x)=V_i(x) = \sum_{\substack{\alpha \in E_{d,m_i} }} b_{i,\alpha} \prod_{j=1}^{m_i} (h^T_{i,j}x+g_{i,j})^{\alpha_j} \; \text{ for } x \in D_i, i=1, \cdots,L
\]
is a piecewise polynomial Lyapunov function proving stability of System~\eqref{eq:sys_nonlin}. See~\cite{reza_CDC_2014} for a comprehensive discussion on the computational complexity of the LP defined in~\eqref{eq:lin_prog} in terms of the state-space dimension and the degree of $V(x)$.

\begin{figure}[htbp]
\includegraphics[scale=0.43]{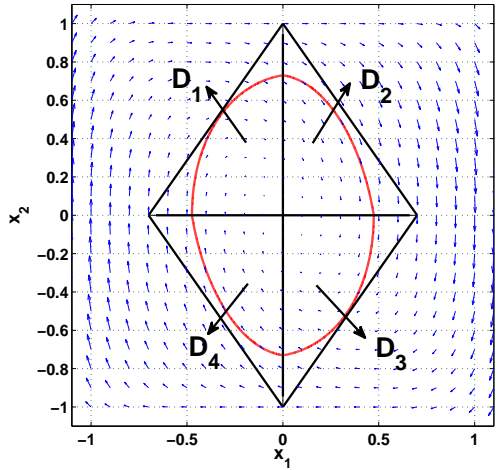}
\caption{The largest level-set of Lyapunov function~\eqref{eq:handelman_PW_Lyap} inscribed in Polytope~\eqref{eq:polytope0}}
\label{fig:handelman_nonlin0}
\end{figure}

\noindent \textbf{\\Numerical Example 4:} Consider the following nonlinear system~\cite{chesi_2005}.
\begin{align*}
& \dot{x}_1 = x_2 \\
& \dot{x}_2 = -2x_1-x_2+x_1 x_2^2-x_1^5+x_1 x_2^4+x_2^5.
\end{align*}
Using the polytope
\begin{align}
\Gamma^4 = \{  x_1,x_2 \in \mathbb{R}^2: \, &1.428 x_1 + x_2 - 0.625 \geq 0, -1.428x_1+x_2+0.625 \geq 0, \nonumber \\
              & 1.428 x_1 + x_2 + 0.625 \geq 0, -1.428x_1+x_2 - 0.625 \geq0 \},
\label{eq:polytope0}
\end{align}
and $D-$decomposition
\begin{align*}
&D_1:=\{x_1,x_2 \in \mathbb{R}^2: -x_1 \geq 0, x_2 \geq 0, -1.428x_1+x_2 - 0.625 \geq 0 \} \\
&D_2:=\{x_1,x_2 \in \mathbb{R}^2: x_1 \geq 0, x_2 \geq 0, 1.428 x_1 + x_2 + 0.625 \geq 0 \}\\
&D_3:=\{x_1,x_2 \in \mathbb{R}^2: x_1 \geq 0, -x_2 \geq 0, -1.428x_1+x_2+0.625 \geq 0  \}\\
&D_4:=\{x_1,x_2 \in \mathbb{R}^2: -x_1 \geq 0, -x_2 \geq 0, 1.428 x_1 + x_2 + 0.625 \geq 0  \},
\end{align*}
we set-up the LP in~\eqref{eq:lin_prog} with $d=4$. The solution to the LP certified asymptotic stability of the origin and yielded the following piecewise polynomial Lyapunov function. Figure~\ref{fig:handelman_nonlin0} shows the largest level set of $V(x)$ inscribed in the polytope $\Gamma^4$. \vspace{-0.4in}
\begin{small}
\begin{align}
V(x) = 
\begin{cases} 
0.543 x_1^2  + 0.233 x_2^2 + 0.018 x_2^3 - 0.074 x_1 x_2^2  - 0.31 x_1^3 \\
  + 0.004 x_2^4 - 0.013 x_1 x_2^3 + 0.015 x_1^2 x_2^2 + 0.315 x_1^4  &\text{if } x \in D_1 \vspace*{0.1in} \\
0.543 x_1^2+ 0.329 x_1 x_2 + 0.233 x_2^2 + 0.018 x_2^3 +0.031 x_1 x_2^2 \\
 + 0.086 x_1^2 x_2 +  0.3 x_1^3 + 0.004 x_2^4 +  0.009 x_1 x_2^3 + 0.015 x_1^2 x_2^2\\
  + 0.008 x_1^3 x_2 + 0.315 x_1^4  &\text{if } x \in D_2  \vspace*{0.1in} \\
0.0543 x_1^2 + 0.0233 x_2^2  -0.0018 x_2^3 + 0.0074 x_1 x_2^2 + 0.03 x_1^3 \\
 + 0.004 x_2^4 -0.013 x_1 x_2^3 + 0.015 x_1^2 x_2^2 + 0.315 x_1^4 &\text{if } x\in D_3  \vspace*{0.1in} \\ 
0.543 x_1^2+ 0.329 x_1 x_2 + 0.233 x_2^2 - 0.018 x_2^3 - 0.031 x_1 x_2^2 \\
 - 0.086 x_1^2 x_2 - 0.3 x_1^3 + 0.004 x_2^4 +  0.009 x_1 x_2^3 + 0.015 x_1^2 x_2^2 \\
  +  0.008 x_1^3 x_2 + 0.315 x_1^4 &\text{if } x \in D_4  
  \label{eq:handelman_PW_Lyap}
\end{cases}
\end{align}
\end{small}
\hspace{-0.12in} See Numerical Example 5 for a comparison of this method with the method in Section~\ref{sec:polya_nonlin}.

\subsubsection{Application of Polya's theorem in nonlinear stability analysis}
\label{sec:polya_nonlin}
In this section, we discuss an algorithm based on a multi-simplex version of Polya's theorem (Theorem~\ref{thm:polya_multi-simplex}) to verify local stability of nonlinear systems of the form
\begin{equation}
\dot{x} = A(x) x(t),
\label{eq:sys_nonlin2}
\end{equation}
where $A(x) \in \mathbb{R}^{n \times n}$ is a matrix-valued polynomial and $A(0) \neq 0$. 

Unfortunately, Polya's theorem does not parameterize polynomials which have zeros in the interior of the unit simplex (see~\cite{polya_corner} for an elementary proof of this). From the same reasoning as in~\cite{polya_corner} it follows that the multi-simplex version of Polya's theorem (Theorem~\ref{thm:polya_multi-simplex}) does not parameterize polynomials which have zeros in the interior of a multi-simplex. On the other hand, if $f(z)$ in~\eqref{eq:fz_multi} has a zero in the interior of $\Phi^n$, then any multi-homogeneous polynomial $p(x,y)$ that satisfies~\eqref{eq:fz_multi} has a zero in the interior of the multi-simplex $\Delta^2 \times \cdots \times \Delta^2$. One way to enforce the condition $V(0) = 0$ in~\eqref{eq:V} is to search for coefficients of a matrix-valued polynomial $P$ which defines a Lyapunov function of the form $V(x)=x^T P(x) x$. It can be shown that $V(x)=x^T P(x) x$ is a Lyapunov function for System~\eqref{eq:sys_nonlin2} if and only if $\gamma^*$ in the following optimization of polynomials problem is positive.
\begin{align}
&\gamma^* =  \max_{\gamma \in \mathbb{R}, P \in \mathbb{R}[x]} \;\; \gamma \nonumber \\
& \text{subject to} 
\begin{bmatrix}
P(x) & 0 \\ 
0 & -Q(x)
\end{bmatrix} - \gamma I \geq 0  \;\; \text{ for all } x \in \Phi^n,
\label{eq:OOP_nonlin2}
\end{align}
where 
\[
Q(x)=
 A^T(x) P(x) + P(x) A(x) + \frac{1}{2} \left( \hspace*{-0.04in}  A^T(x) \begin{bmatrix}
x^T \frac{\partial P(x)}{\partial x_1} \\ 
\vdots \\ 
x^T \frac{\partial P(x)}{\partial x_n}
\end{bmatrix}
+
 \begin{bmatrix}
x^T \frac{\partial P(x)}{\partial x_1} \\ 
\vdots \\ 
x^T \frac{\partial P(x)}{\partial x_n}
\end{bmatrix}^T A(x)  \right).
\]
As in Case 2 of Section~\ref{sec:algorithms}, by applying bisection algorithm on $\gamma$ and using Theorem~\ref{thm:polya_multi-simplex} as a test for feasibility of Constraint~\eqref{eq:OOP_nonlin2} we can compute lower bounds on $\gamma^*$. Suppose there exists a matrix-valued multi-homogeneous  polynomial $S$ of degree vector $d_s \in \mathbb{N}^n$ ($d_{s_i}$ is the degree of $y_i$) such that
\[
\{ P(x)\in \mathbb{S}^n: x \in \Phi^n \} = \{ S(y,z) \in \mathbb{S}^n: (y_i,z_i) \in \Delta^2, i=1,\cdots,n \}.
\]
Likewise, suppose there exist matrix-valued multi-homogeneous polynomials $B$ and $C$ of degree vectors $d_b \in \mathbb{N}^n$ and $d_c=d_s \in \mathbb{N}^n$ such that
\[
\{ A(x) \in \mathbb{R}^{n \times n}: x \in \Phi^n \} = \{ B(y,z) \in \mathbb{R}^{n \times n}: (y_i,z_i) \in \Delta^2, i=1,\cdots,n \}
\]
and 
\begin{small}
\[
\left\lbrace  \begin{bmatrix}
 \frac{\partial P(x)}{\partial x_1} x, \cdots , \frac{\partial P(x)}{\partial x_n}  x
\end{bmatrix} \in \mathbb{R}^{n \times n}: x \in \Phi^n \right\rbrace \hspace*{-0.03in}  =  \hspace*{-0.03in} \left\lbrace C(y,z) \in \mathbb{R}^{n \times n}: (y_i,z_i) \in \Delta^2, i=1,\cdots,n \right\rbrace \hspace*{-0.02in} .
\]
\end{small}
Given $\gamma \in \mathbb{R}$, it follows from Theorem~\ref{thm:polya_multi-simplex} that the inequality condition in~\eqref{eq:OOP_nonlin2} holds for all $\alpha \in \Phi^l$ if there exist $e \geq 0$ such that 
\begin{equation}
\left( \prod_{i=1}^n (y_i+z_i)^e \right) \left( S(y,z) - \gamma I \left( \prod_{i=1}^n (y_i+z_i)^{d_{p_i}} \right) \right)
\label{eq:product1}
\end{equation}
and
\begin{align}
\left( \prod_{i=1}^n (y_i+z_i)^e \right) \left(
 B^T(y,z) S(y,z) + S(y,z) B(y,z) \right. \nonumber \\
 & \hspace{-2.6in} \left. + \frac{1}{2} \left(   B^T(y,z) C^T(y,z) +
C(y,z) B(y,z)  \right) - \gamma I  \left( \prod_{i=1}^n (y_i+z_i)^{d_{q_i}} \right)  \right)
\label{eq:product2}
\end{align}
have all positive coefficients, where $d_{p_i}$ is the degree of $x_i$ in $P(x)$ and $d_{q_i}$ is the degree of $x_i$ in $Q(x)$. Suppose $S,B$ and $C$ have the following forms. 
\begin{equation}
S(y,z) = \sum_{ \substack{ h,g \in \mathbb{N}^l \\ h+g=d_s }} S_{h,g} y_1^{h_1} z_1^{g_1} \cdots y_n^{h_l} z_n^{g_l},
\label{eq:S_polya_nonlin}
\end{equation}
\begin{equation}
B(y,z) = \sum_{ \substack{ h,g \in \mathbb{N}^l \\ h+g=d_b }} B_{h,g} y_1^{h_1} z_1^{g_1} \cdots y_n^{h_l} z_n^{g_l}
\label{eq:B_polya_nonlin}
\end{equation}
\begin{equation}
C(y,z) = \sum_{ \substack{ h,g \in \mathbb{N}^l \\ h+g=d_c }} C_{h,g} y_1^{h_1} z_1^{g_1} \cdots y_n^{h_l} z_n^{g_l},
\label{eq:C_polya_nonlin}
\end{equation}
By combining~\eqref{eq:S_polya_nonlin},~\eqref{eq:B_polya_nonlin} and~\eqref{eq:C_polya_nonlin} with~\eqref{eq:product1} and~\eqref{eq:product2} it follows that for a given $\gamma \in \mathbb{R}$, the inequality condition in~\eqref{eq:OOP_nonlin2} holds for all $\alpha \in \Phi^n$ if there exist some $e \geq 0$ such that
\begin{align}
& \sum_{\substack{ h,g \in \mathbb{N}^l \\ h+g=d_s }} f_{\{q,r\}, \{h,g \}} S_{h,g} > 0 \;\; \text{ for all } \; q,r \in \mathbb{N}^l: q+r=d_s+e \cdot \mathbf{1} \;\; \text{ and} \nonumber \\
& \sum_{\substack{ h,g \in \mathbb{N}^l \\ h+g=d_s }} \hspace*{-0.09in}  M_{\{u,v\},\{h,g\}}^T S_{h,g} + S_{h,g} M_{\{u,v\},\{h,g\}}  +  N_{\{u,v\},\{h,g\}}^T C^T_{h,g} + C_{h,g} N_{\{u,v\},\{h,g\}} < 0  \nonumber \\
& \hspace*{2.5in} \text{ for all }  u,v \in \mathbb{N}^l: u+v=d_s+d_b+e \cdot \mathbf{1},
\label{eq:polya_nonlin_conditions}
\end{align}
where similar to Case 3 of Section~\ref{sec:applications_lin}, we define $f_{\{ q,r \},\{h,g\}}$ to be the coefficient of $S_{h,g} y^qz^r$ after combining~\eqref{eq:S_polya_nonlin} with~\eqref{eq:product1}. Likewise, we define $M_{\{u,v\},\{h,g\}}$ to be the coefficient of $S_{h,g} y^uz^v$ and $N_{\{u,v\},\{h,g\}}$ to be the coefficient of $C_{h,g} y^uz^v$ after combining~\eqref{eq:B_polya_nonlin} and~\eqref{eq:C_polya_nonlin} with~\eqref{eq:product2}.
Conditions~\eqref{eq:polya_nonlin_conditions} are an SDP (See~\cite{reza_CDC_2014} for a complexity analysis on this SDP). For any $\gamma \in \mathbb{R}$, if there exist $e \geq 0$ and $\{ S_{h,g} \}$ such that Conditions~\eqref{eq:polya_nonlin_conditions} hold, then $\gamma$ is a lower bound for $\gamma^*$ as defined in~\eqref{eq:OOP_nonlin2}. Furthermore, if $\gamma$ is positive, then origin is an asymptotically stable equilibrium for System~\eqref{eq:sys_nonlin2}.

\begin{figure}[htbp]
\includegraphics[scale=0.26]{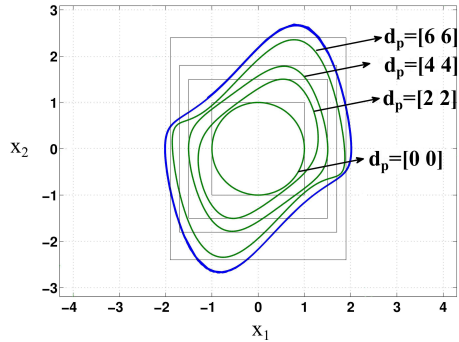}
\caption{Level-sets of the Lyapunov functions $V(x)=x^T P(x) x$ computed by the method in Section~\ref{sec:polya_nonlin} - $d_p$ is the degree of $P(x)$ } 
\label{fig:polya_nonlin}
\end{figure}

\noindent \textbf{\\Numerical Example 5:} Consider the reverse-time Van Der Pol oscillator defined as
\[
\begin{bmatrix}
\dot{x}_1  \\ 
\dot{x}_2
\end{bmatrix} 
=
\begin{bmatrix}
-x_2  \\ 
x_1+x_2(x_1^2-1)
\end{bmatrix} 
=
A(x)x,
\]
where $A(x)=\begin{bmatrix}
0 & -1 \\ 
1 & x_1^2-1
\end{bmatrix} $. By using the method in Section~\ref{sec:polya_nonlin}, we solved Problem~\eqref{eq:OOP_nonlin2} using the hypercubes
\begin{align}
& \Phi_1^2 = \{ x \in \mathbb{R}^2: \vert x_1 \vert \leq 1, \vert x_2 \vert \leq 1\} \nonumber \\
& \Phi_2^2 = \{ x \in \mathbb{R}^2: \vert x_1 \vert \leq 1.5, \vert x_2 \vert \leq 1.5 \} \nonumber \\
& \Phi_3^2 = \{ x \in \mathbb{R}^2: \vert x_1 \vert \leq 1.7, \vert x_2 \vert \leq 1.8 \} \nonumber \\
& \Phi_4^2 = \{ x \in \mathbb{R}^2: \vert x_1 \vert \leq 1.9, \vert x_2 \vert \leq 2.4\}
\label{eq:hypercubes}
\end{align}
and $d_p=0,2,4,6$ as the degrees of $P(x)$. For each hypercube $\Phi_i^2$ in~\eqref{eq:hypercubes}, we computed a Lyapunov function of the form $V_i(x)=x^T P_i(x) x$. In Figure~\ref{fig:polya_nonlin}, we have plotted the largest level-set of $V_i$, inscribed in $\Phi_i^2$ for $i=1, \cdots,4$. For all the cases, we used the Polya's exponent $e=1$.


\begin{figure}[htbp]
\includegraphics[scale=0.37]{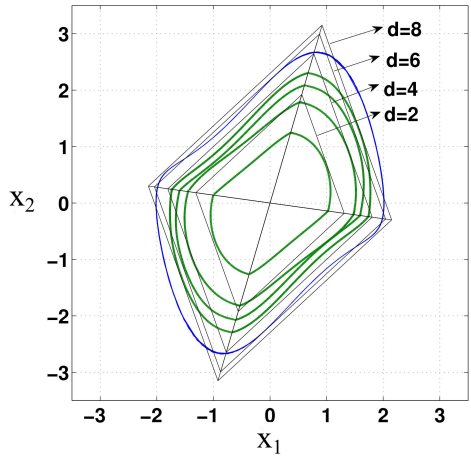}
\caption{Level-sets of the Lyapunov functions computed for Van Der Pol oscillator using the method in Section~\ref{sec:handelman_nonlin} - $d$ is the degree of the Lyapunov functions}
\label{fig:handelman_nonlin} 
\end{figure}

We also used the method in Section~\ref{sec:handelman_nonlin} to solve  the same problem (see Figure~\ref{fig:handelman_nonlin}) using the polytopes
\[
\Gamma_\nu:= \left\lbrace  x \in \mathbb{R}^2: x= \sum_{i=1}^4 \rho_i v_i: \rho_i \in [0,\nu], \sum_{i=1}^4\rho_i = \nu  \right\rbrace 
\]
with $\nu=0.83, 1.41 ,1.52,1.64$, where
\[
v_1= \begin{bmatrix}
-1.31 \\
0.18
\end{bmatrix},
v_2=\begin{bmatrix}
0.56\\
1.92
\end{bmatrix},
v_3=\begin{bmatrix}
-0.56\\
-1.92
\end{bmatrix} \text{ and }
v_4=\begin{bmatrix}
1.31\\
-0.18
\end{bmatrix}.
\]

From Figures~\ref{fig:polya_nonlin} and~\ref{fig:handelman_nonlin} we observe that in both methods, computing larger invariant subsets of the region of attraction of the origin requires an increase in the degree of Lyapunov functions. The CPU time required for computing the Lyapunov functions associated with the largest invariant subsets in Figures~\ref{fig:polya_nonlin} and~\ref{fig:handelman_nonlin} were 88.9 minutes (using the method of Section~\ref{sec:polya_nonlin}) and 3.78 minutes (using the method of Section~\ref{sec:handelman_nonlin}), respectively using a core i7 machine with 64 GB of RAM.

\subsection{Robust $H_{\infty}$ control synthesis}

Consider plant $G$ with the state-space formulation 
\begin{equation}
\dot{x}(t) = A(\alpha) x(t) + \begin{bmatrix}
B_1(\alpha) & B_2(\alpha)
\end{bmatrix}
\begin{bmatrix}
\omega(t) \nonumber \\
u(t)
\end{bmatrix}
\end{equation}
\begin{equation}
\begin{bmatrix}
z(t) \\
y(t)
\end{bmatrix} =
\begin{bmatrix}
C_1(\alpha) \\
C_2(\alpha)
\end{bmatrix} x(t) +
\begin{bmatrix}
D_{11}(\alpha) & D_{12}(\alpha) \\
D_{21}(\alpha) & 0
\end{bmatrix}
\begin{bmatrix}
\omega(t) \\
u(t)
\end{bmatrix},
\label{eq:sys_G}
\end{equation}
where $\alpha \in Q \subset \mathbb{R}^l$, $x(t) \in \mathbb{R}^n$, $u(t) \in \mathbb{R}^m$, $\omega(t) \in \mathbb{R}^p$ is the external input and $z(t) \in \mathbb{R}^q$ is the external output. Suppose $(A(\alpha),B_2(\alpha))$ is stabilizable and $(C_2(\alpha),A(\alpha))$ is detectable for all $\alpha \in Q$.
 According to~\cite{Gahinet_1994_LMI_hinf} there exists a state feedback gain $K(\alpha) \in \mathbb{R}^{m \times n}$ such that
\[
\| S (G,K(\alpha))  \|_{H_{\infty}} \leq \gamma, \; \text{for all} \; \alpha \in Q,
\]
if and only if there exist $P(\alpha) > 0$ and $R(\alpha) \in \mathbb{R}^{m \times n}$ such that $K(\alpha)=R(\alpha)P^{-1}(\alpha)$ and
\begin{small}
\begin{equation}
\begin{bmatrix}

\begin{bmatrix}
A(\alpha) \hspace*{-0.05in} &  B_2(\alpha)
\end{bmatrix} \hspace*{-0.05in} \begin{bmatrix}
P(\alpha) \\
R(\alpha)
\end{bmatrix}  \hspace*{-0.05in} + \hspace*{-0.05in} \begin{bmatrix}
P(\alpha) \hspace*{-0.05in} & R^T(\alpha)
\end{bmatrix} \hspace*{-0.05in} \begin{bmatrix}
A^T(\alpha) \\
B_2^T(\alpha)
\end{bmatrix} & \star & \star \\
B_1^T(\alpha) & -\gamma I & \star \\
\begin{bmatrix}
C_1(\alpha) & D_{12}(\alpha)
\end{bmatrix}\begin{bmatrix}
P(\alpha) \\
R(\alpha)
\end{bmatrix} & D_{11}(\alpha) & -\gamma I
\end{bmatrix} < 0,
\end{equation}
\end{small}
for all $\alpha \in Q$, where $\gamma > 0$ and $S(G,K(\alpha))$ is the map from the external input $\omega$ to the external output $z$ of the closed loop system with a static full state feedback controller. The symbol $\star$ denotes the symmetric blocks in the matrix inequality. 
To find a solution to the robust $H_\infty$-optimal static state-feedback controller problem with optimal feedback gain $K(\alpha)=P(\alpha)R^{-1}(\alpha)$, one can solve the following optimization of polynomials problem.
\begin{small}
\begin{align}
&\gamma^* = \min_{P,R \in \mathbb{R}[\alpha],\gamma \in \mathbb{R}} \;\; \gamma \nonumber \\
& \text{subject to } \nonumber \\
& \hspace*{-0.07in}\begin{bmatrix}
-P(\alpha) & \star  & \hspace*{-0.1in} \star  & \star \\
0 & \hspace*{-0.05in}
\begin{bmatrix}
 A(\alpha) \hspace*{-0.09in} &  B_2(\alpha)
\end{bmatrix} \hspace*{-0.05in} \begin{bmatrix}
P(\alpha) \\
R(\alpha)
\end{bmatrix}  \hspace*{-0.05in} + \hspace*{-0.05in} \begin{bmatrix}
P(\alpha) \hspace*{-0.05in} & R^T(\alpha)
\end{bmatrix} \hspace*{-0.05in} \begin{bmatrix}
A^T(\alpha) \\
B_2^T(\alpha)
\end{bmatrix} & \hspace*{-0.1in} \star &  \star\\
0 & 
B_1^T(\alpha) & \hspace*{-0.1in} 0 & \star \\
0 & 
\begin{bmatrix}
C_1(\alpha) & D_{12}(\alpha)
\end{bmatrix}
\begin{bmatrix}
P(\alpha) \\
R(\alpha)
\end{bmatrix} & \hspace*{-0.1in} D_{11}(\alpha) & 0
\end{bmatrix} 
& \hspace*{-0.15in} -\gamma 
\begin{bmatrix}
0 & 0 & 0 & 0 \\ 
0 & 0 & 0 & 0 \\ 
0 & 0 & I & 0 \\ 
0 & 0 & 0 & I
\end{bmatrix} < 0 \nonumber \\
& \hspace*{3in} \text{ for all }  \alpha \in Q.
\label{eq:hinf_ineq}
\end{align}
\end{small}

In Problem~\eqref{eq:hinf_ineq}, if $Q = \Delta^l$ as defined in~\eqref{eq:simplex}, then we can apply Polya's theorem (Theorem~\ref{thm:polya_simplex}) as in the algorithm in Case 1 of Section~\ref{sec:algorithms} to find a $\gamma \leq \gamma^*$ and $P$ and $R$ which satisfy the inequality in~\eqref{eq:hinf_ineq}. Suppose $P, A, B_1, B_2, C_1, D_{11}$ and $D_{12}$ are homogeneous polynomials. If any of these polynomials is not homogeneous, use the procedure in Case 1 of Section~\ref{sec:algorithms} to homogenize it. Let
\begin{small}
\[
F(P(\alpha), R(\alpha)) := \begin{bmatrix}
-P(\alpha) & \star  & \hspace*{-0.1in} \star  & \star \\
0 & \hspace*{-0.05in}
\begin{bmatrix}
 A(\alpha) \hspace*{-0.09in} &  B_2(\alpha)
\end{bmatrix} \hspace*{-0.05in} \begin{bmatrix}
P(\alpha) \\
R(\alpha)
\end{bmatrix}  \hspace*{-0.05in} + \hspace*{-0.05in} \begin{bmatrix}
P(\alpha) \hspace*{-0.05in} & R^T(\alpha)
\end{bmatrix} \hspace*{-0.05in} \begin{bmatrix}
A^T(\alpha) \\
B_2^T(\alpha)
\end{bmatrix} & \hspace*{-0.1in} \star &  \star\\
0 & 
B_1^T(\alpha) & \hspace*{-0.1in} 0 & \star \\
0 & 
\begin{bmatrix}
C_1(\alpha) & D_{12}(\alpha)
\end{bmatrix}
\begin{bmatrix}
P(\alpha) \\
R(\alpha)
\end{bmatrix} & \hspace*{-0.1in} D_{11}(\alpha) & 0
\end{bmatrix},
\]
\end{small}
\hspace*{-0.08in} and denote the degree of $F$ by $d_f$. Given $\gamma \in \mathbb{R}$, the inequality in~\eqref{eq:hinf_ineq} holds if there exist $e \geq 0$ such that all of the coefficients in
\begin{equation}
\left( \sum_{i=1}^l \alpha_i \right)^e \left( F(P(\alpha), R(\alpha)) - \gamma 
\begin{bmatrix}
0 & 0 & 0 & 0 \\ 
0 & 0 & 0 & 0 \\ 
0 & 0 & I & 0 \\ 
0 & 0 & 0 & I
\end{bmatrix} \left( \sum_{i=1}^l \alpha_i \right)^{d_f} \right)
\label{eq:polya_hinf}
\end{equation}
are negative-definite. Let $P$ and $R$ be of the forms
\begin{equation}
P(\alpha) = \sum_{h \in I_{d_p}} P_h \alpha_1^{h_1} \cdots \alpha_l^{h_l}, P_h \in \mathbb{S}^n \text{ and } R(\alpha) = \sum_{h \in I_{d_r}} R_h \alpha_1^{h_1} \cdots \alpha_l^{h_l}, R_h \in \mathbb{R}^{n \times n},
\label{eq:PR}
\end{equation}
where $I_{d_p}$ and $I_{d_r}$ are defined as in~\eqref{eq:index_set}. By combining~\eqref{eq:PR} with~\eqref{eq:polya_hinf} it follows from Theorem~\eqref{thm:polya_simplex} that for a given $\gamma$, the inequality in~\eqref{eq:hinf_ineq} holds, if there exist $e \geq 0$ such that
\begin{equation}
\sum_{h \in I_{d_p}}\left( M_{h,q}^T P_h+P_h M_{h,q} \right) + \sum_{h \in I_{d_r}}\left( N_{h,q}^T R_h^T + R_h N_{h,q} \right) < 0 \;\; \text{ for all }  q \in I_{d_f+e},
\label{eq:LMI_hinf}
\end{equation}
where we define $M_{h,q} \in \mathbb{R}^{n \times n}$ as the coefficient of $P_h \alpha_1^{q_1} \cdots \alpha_l^{q_l}$ after combining~\eqref{eq:PR} with~\eqref{eq:polya_hinf}. Likewise, $N_{h,q} \in \mathbb{R}^{n \times n}$ is the coefficient of $R_h \alpha_1^{q_1} \cdots \alpha_l^{q_l}$ after combining~\eqref{eq:PR} with~\eqref{eq:polya_hinf}. For given $\gamma > 0$, if there exist $e \geq 0$ such that the LMI~\eqref{eq:LMI_hinf} has a solution, say $P_h, h\in I_{d_p}$ and $R_g, g\in I_{d_r}$, then 
\[
K(\alpha)= \left( \sum_{h \in I_{d_p}}  P_h \alpha_1^{h_1} \cdots \alpha_l^{h_l} \right) \left( \sum_{g \in I_{d_r}}  R_g \alpha_1^{g_1} \cdots \alpha_l^{g_l} \right)^{-1}
\]
is a feedback gain for an $H_{\infty}$-suboptimal static state-feedback controller 
for System~\eqref{eq:sys_G}. By performing a bisection on $\gamma$ and solving~\eqref{eq:LMI_hinf} form each $\gamma$ of the bisection, one may find an $H_{\infty}$-optimal controller for System~\eqref{eq:sys_G}.

In Problem~\eqref{eq:hinf_ineq}, if $Q = \Phi^l$ as defined in~\eqref{eq:hypercube}, then by applying the algorithm in Case 2 of section~\ref{sec:algorithms} to Problem~\eqref{eq:hinf_ineq}, we can find a solution $P, Q, \gamma$ to~\eqref{eq:hinf_ineq}, where $\gamma  \leq \gamma^*$. See Case 3 of Section~\ref{sec:applications_lin} and Section~\ref{sec:polya_nonlin} for similar applications of this theorem.

If $Q= \Gamma^l$ as defined in~\eqref{eq:polytope}, then we can use Handelman's theorem (Theorem~\ref{thm:Handelman}) as in the algorithm in Case 3 of section~\ref{sec:algorithms} to find a solution to Problem~\eqref{eq:hinf_ineq}. We have provided a similar application of Handelman's theorem in Section~\ref{sec:handelman_nonlin}.

If $Q$ is a compact semi-algebraic set, then for given $d \in \mathbb{N}$, one can apply the Positivstellensatz results in Case 4 of Section~\ref{sec:algorithms} to the inequality in~\eqref{eq:hinf_ineq} to obtain a SOS program of the Form~\eqref{eq:SOS_putinar}. A solution to the SOS program yields a solution to Problem~\eqref{eq:hinf_ineq}.


\section{Conclusion}
SOS programming, moment's approach and their applications in polynomial optimization have been well-served in the literature. To promote diversity in commonly used algorithms for polynomial optimization, we dedicated this paper to some of the alternatives to SOS programming. In particular, we focused on the algorithms defined by Polya's theorem, Bernstein's theorem and Handelman's theorem. We discussed how these algorithms apply to polynomial optimization problems with decision variables on  simplices, hypercubes and arbitrary convex polytopes. Moreover, we demonstrated some of the applications of Polya's and Handelman's algorithms in stability analysis of nonlinear systems and stability analysis and $H_{\infty}$ control of systems with parametric uncertainty. For most of these applications, we have provided numerical examples to compare the conservativeness of Polya's and Handelman's algorithms with other algorithms in the literature such as SOS programming.\\



\textbf{Acknowledgements.} This material is based upon work supported by the National Science Foundation under Grant Number 1301660.


\begin{thebibliography}{99}

\bibitem{blum} (MR1301779) [10.1007/978-1-4612-0873-0]
     \newblock L. Blum, F. Cucker, M. Shub and S. Smale,
     \newblock \emph{Complexity and real computation},
     \newblock Springer-Verlag, New York, 1998.
     
\bibitem{inequalities} (MR1301779) [10.1007/978-1-4612-0873-0]
     \newblock G. Hardy, J. Littlewood and G. Polya,
     \newblock \emph{Inequalities},
     \newblock Cambridge University Press, 1934.
     
\bibitem{adams_groebner} (MR1301779) [10.1007/978-1-4612-0873-0]
     \newblock W. Adams and P. Loustaunau,
     \newblock \emph{An Introduction to Groebner Bases},
     \newblock American Mathematical Society, 1994.
     
\bibitem{blossoming} (MR1301779) [10.1007/978-1-4612-0873-0]
     \newblock L. Ramshaw,
     \newblock \emph{A connect-the-dots approach to splines},
     \newblock Digital Systems Research Center, 1987.

\bibitem{tarski} (MR1124979) [10.2307/2152750]
     \newblock A. Tarski,
     \newblock \emph{A Decision Method for Elementary Algebra and Geometry},
     \newblock Random Corporation monograph, Berekley and Los Angeles, 1951.      
     
\bibitem{sherali_global_2007} (MR1124979) [10.2307/2152750]
    \newblock H. Sherali and L. Liberti,
    \newblock \emph{Reformulation-linearization technique for global optimization},
    \newblock Encyclopedia of Optimization, Springer, USA, 2009, 3263--3268. 

\bibitem{parillo_thesis}
    \newblock P. Parrilo,
    \newblock  \emph{Structured semidefinite programs and semialgebraic geometry methods in robustness and optimization},
    \newblock  Ph.D thesis, California Institute of Technology, 2000.

\bibitem{sostools2013}
\newblock A. Papachristodoulou, J. Anderson, G. Valmorbida, S. Prajna, P. Seiler and P. A. Parrilo,
\newblock SOSTOOLS: Sum of squares optimization toolbox for MATLAB,
\newblock preprint, \arXiv{1310.4716}, 2013.

\bibitem{monteiro} (MR1124979) [10.2307/2152750]
    \newblock R. Monteiro
    \newblock Primal-Dual Path-Following Algorithms for Semidefinite Programming,
    \newblock \emph{SIAM Journal of Optimization}, \textbf{7(3)} (1997), 663--678.
    
\bibitem{helmberg} (MR1124979) [10.2307/2152750]
    \newblock C. Helmberg, F. Rendl, R. J. Vanderbei and H. Wolkowicz,
    \newblock An Interior-Point Method for Semidefinite Programming,
    \newblock \emph{SIAM Journal of Optimization}, \textbf{6(2)} (1996), 342--361.
    
\bibitem{alizadeh} (MR1124979) [10.2307/2152750]
    \newblock F. Alizadeh, J. Haeberly and M. Overton,
    \newblock Primal-Dual Interior-Point Methods for Semidefinite Programming: Convergence Rates, Stability and Numerical Results,
    \newblock \emph{SIAM Journal of Optimization}, \textbf{8(3)} (1998), 746--768.

\bibitem{sdpa} (MR1124979) [10.2307/2152750]
    \newblock M. Yamashita et. al.,
    \newblock A high-performance software package for semidefinite programs: SDPA 7,
    \newblock \emph{Tech. rep. B-460, Dep. of Mathematical and Computing Sciences, Tokyo Inst. of Tech.}, (2010).

\bibitem{sdpt3} (MR1124979) [10.2307/2152750]
    \newblock R. Tutuncu, K. Toh and M. Todd,
    \newblock Solving semidefinite-quadratic-linear programs using SDPT3, Mathematical Programming,
    \newblock \emph{Mathematical Programming Series B}, \textbf{95} (2003), 189--217.

\bibitem{sedumi} (MR1124979) [10.2307/2152750]
    \newblock J. Sturm,
    \newblock Using SeDuMi 1.02, a MATLAB toolbox for optimization over symmetric cones,
    \newblock \emph{Optimization methods and software}, \textbf{11(1-4)} (1999), 625--653.
    
\bibitem{stengle} (MR1124979) [10.2307/2152750]
    \newblock G. Stengle,
    \newblock A Nullstellensatz and a Positivstellensatz in semialgebraic geometry,
    \newblock \emph{Mathematische Annalen}, \textbf{207(2)} (1974), 87--97.
    
\bibitem{putinar} (MR1124979) [10.2307/2152750]
    \newblock M. Putinar,
    \newblock Positive polynomials on compact semi-algebraic sets,
    \newblock \emph{Indiana University Mathematics Journal}, \textbf{42} (1993), 969--984.
    
\bibitem{schmudgen} (MR1124979) [10.2307/2152750]
    \newblock K. Schmudgen,
    \newblock The K-moment problem for compact semi-algebraic sets,
    \newblock \emph{Mathematische Annalen}, \textbf{289} (1991), 203--206. 
    
\bibitem{laurent} (MR1124979) [10.2307/2152750]
    \newblock M. Laurent,
    \newblock Sums of squares, moment matrices and optimization over polynomials,
    \newblock \emph{Emerging applications of algebraic geometry}, Springer New York \textbf{11(1-4)} (2009), 157--270.   
      
\bibitem{CAD} (MR1124979) [10.2307/2152750]
    \newblock G. Collins and H. Hoon,
    \newblock Partial cylindrical algebraic decomposition for quantifier elimination,
    \newblock \emph{Journal of Symbolic Computation}, \textbf{12(3)} (1991), 299--328. 
    
\bibitem{sherali_1992} (MR1124979) [10.2307/2152750]
    \newblock H. Sherali and C. Tuncbilek,
    \newblock A global optimization algorithm for polynomial programming problems using a reformulation-linearization technique,
    \newblock \emph{Journal of Global Optimization}, \textbf{2} (1992), 101--112. 
 
\bibitem{sherali_1997} (MR1124979) [10.2307/2152750]
    \newblock H. Sherali and C. Tuncbilek,
    \newblock New reformulation- linearization technique based relaxations
for univariate and multivariate polynomial programming problems,
    \newblock \emph{Operations Research Letters}, \textbf{21(1)} (1997), 1--10.  
    
\bibitem{roy} (MR1124979) [10.2307/2152750]
    \newblock F. Boudaoud, F. Caruso and M Roy,
    \newblock Certificates of Positivity in the Bernstein Basis,
    \newblock \emph{Discrete and Computational Geometry}, \textbf{39(4)} (2008), 639--655.  
    
\bibitem{leroy} (MR1124979) [10.2307/2152750]
    \newblock R. Leroy,
    \newblock Convergence under Subdivision and Complexity of Polynomial Minimization in the Simplicial Bernstein Basis,
    \newblock \emph{Reliable Computing}, \textbf{17} (2012), 11--21. 
    
\bibitem{handelman_1988} (MR1124979) [10.2307/2152750]
    \newblock D. Handelman,
    \newblock Representing polynomials by positive linear functions on compact convex polyhedra,
    \newblock \emph{Pacific Journal of Mathematics}, \textbf{132(1)} (1988), 35--62. 
    
\bibitem{QEPCAD} (MR1124979) [10.2307/2152750]
    \newblock C. Brown,
    \newblock QEPCAD B: a program for computing with semi-algebraic sets using CADs,
    \newblock \emph{ACM SIGSAM Bulletin}, \textbf{37(4)} (2003), 97--108. 

\bibitem{Redlog} (MR1124979) [10.2307/2152750]
    \newblock A. Dolzmann and T. Sturm,
    \newblock Redlog: Computer algebra meets computer logic,
    \newblock \emph{ACM SIGSAM Bulletin}, \textbf{31(2)} (1997), 2--9.

\bibitem{sherali_1990} (MR1124979) [10.2307/2152750]
    \newblock H. Sherali and W. Adams,
    \newblock A hierarchy of relaxations between the continuous and convex hull representations for zero-one programming problems,
    \newblock \emph{SIAM Journal on Discrete Mathematics}, \textbf{3(3)} (1990), 411--430.

\bibitem{PP_groebner} (MR1124979) [10.2307/2152750]
    \newblock Y. Chang and B. Wah,
    \newblock Polynomial Programming Using Groebner Bases,
    \newblock \emph{IEEE Computer Software and Applications Conference}, \textbf{3(3)} (1994), 236--241.

\bibitem{reza_peet_tac2013} (MR1124979) [10.2307/2152750]
    \newblock R. Kamyar, M. Peet, Y. Peet,
    \newblock Solving large-scale robust stability problems by exploiting the parallel structure of Polya's theorem,
    \newblock \emph{IEEE Transactions on Automatic Control}, \textbf{58(8)} (2013), 1931--1947.

\bibitem{girard_blossom} (MR1124979) [10.2307/2152750]
    \newblock M. Ben Sassi and A. Girard,
    \newblock Computation of polytopic invariants for polynomial dynamical systems using linear programming,
    \newblock \emph{Automatica}, \textbf{48(12)} (2012), 3114--3121.
    
\bibitem{polya_corner} (MR1124979) [10.2307/2152750]
    \newblock V. Powers and B. Reznick,
    \newblock A quantitative Polya's Theorem with corner zeros,
    \newblock \emph{ACM International Symposium on Symbolic and Algebraic Computation},  (2006).
    
\bibitem{polya_edge} (MR1124979) [10.2307/2152750]
    \newblock M. Castle, V. Powers and B. Reznick,
    \newblock Polya's theorem with zeros,
    \newblock \emph{Journal of Symbolic Computation}, \textbf{46(9)} (2011), 1039--1048.    

\bibitem{peres_multisimplex} (MR1124979) [10.2307/2152750]
    \newblock R. Oliveira, P. Bliman and P. Peres,
    \newblock Robust LMIs with parameters in multi-simplex: Existence of solutions and applications,
    \newblock \emph{IEEE 47th Conference on Decision and Control}, (2008), 2226--2231.
    
\bibitem{reza_CDC_hypercube} (MR1124979) [10.2307/2152750]
    \newblock R. Kamyar and M. Peet,
    \newblock Decentralized computation for robust stability of large-scale systems with parameters on the hypercube,
    \newblock \emph{IEEE 51st Conference on Decision and Control}, (2012), 6259--6264.          
    
\bibitem{polya_positivstellensatz_2014} (MR1124979) [10.2307/2152750]
    \newblock P. Dickinson and J. Pohv,
    \newblock On an extension of Polya's Positivstellensatz,
    \newblock \emph{Journal of Global Optimization} (2014), 1--11.      
    
\bibitem{polya_Rn} (MR1124979) [10.2307/2152750]
    \newblock J. de Loera and F. Santos,
    \newblock An effective version of Polya's theorem on positive definite forms,
    \newblock \emph{Journal of Pure and Applied Algebra}, \textbf{108(3)} (1996), 231--240.    
    
\bibitem{polya_rational} (MR1124979) [10.2307/2152750]
    \newblock C. Delzell,
    \newblock Impossibility of extending Polya's theorem to forms with arbitrary real exponents,
    \newblock \emph{Journal of Pure and Applied Algebra}, \textbf{212(12)} (2008), 2612--2622.
    
\bibitem{peres_ramos_2002} (MR1124979) [10.2307/2152750]
    \newblock D. Ramos and P. Peres,
    \newblock An LMI Condition for the Robust Stability of Uncertain
Continuous-Time Linear Systems,
    \newblock \emph{IEEE Transactions on Automatic Control}, \textbf{47(4)} (2002), 675--678.   
    
\bibitem{peres_TAC2007} (MR1124979) [10.2307/2152750]
    \newblock R. Oliveira and P. Peres,
    \newblock Parameter-dependent LMIs in robust analysis: characterization of homogeneous polynomially parameter-dependent solutions via LMI relaxations,
    \newblock \emph{IEEE Transactions on Automatic Control}, \textbf{52(7)} (2007), 1334--1340.
    
\bibitem{chesi_ppdlf_2005} (MR1124979) [10.2307/2152750]
    \newblock G. Chesi, A. Garulli, A. Tesi and A. Vicino,
    \newblock Polynomially parameter-dependent Lyapunov functions for robust stability of polytopic systems: an LMI approach,
    \newblock \emph{IEEE Transactions on Automatic Control}, \textbf{50(3)} (2005), 365--370.    
    
\bibitem{CPA_sigurdur} (MR1124979) [10.2307/2152750]
    \newblock P. Giesel and H. Sigurdur,
    \newblock Revised CPA method to compute Lyapunov functions for nonlinear systems,
    \newblock \emph{Journal of Mathematical Analysis and Applications}, \textbf{410(1)} (2014), 292--306.     
    
\bibitem{reza_CDC_2013} (MR1124979) [10.2307/2152750]
    \newblock R. Kamyar and M. Peet,
    \newblock Decentralized polya's algorithm for stability analysis of large-scale nonlinear systems,
    \newblock \emph{IEEE Conference on Decision and Control}, (2013), 5858--5863.
    

\bibitem{Handelman_Sankaranarayanan} (MR1124979) [10.2307/2152750]
    \newblock S. Sankaranarayanan, X. Chen and E. Abrahám,
    \newblock Lyapunov Function Synthesis using Handelman Representations,
    \newblock \emph{The 9th IFAC Symposium on Nonlinear Control Systems}, (2013).   
    
\bibitem{reza_CDC_2014} (MR1124979) [10.2307/2152750]
    \newblock R. Kamyar, C. Murti and M. Peet,
    \newblock Constructing Piecewise-Polynomial Lyapunov Functions on Convex
Polytopes Using Handelman's Basis,
    \newblock \emph{IEEE Conference on Decision and Controls}, (2014).   

    
\bibitem{peet_SIAM_delay2009} (MR1124979) [10.2307/2152750]
    \newblock A. Papachristodoulou, M. Peet and S. Lall,
    \newblock Analysis of polynomial systems with time delays via the sum of squares decomposition,
    \newblock \emph{IEEE Transactions on Automatic Control}, \textbf{54(5)} (2009), 1058--1064.     
    
\bibitem{peet_TAC_delay2009} (MR1124979) [10.2307/2152750]
    \newblock M. Peet, A. Papachristodoulou and S. Lall,
    \newblock Positive forms and stability of linear time-delay systems,
    \newblock \emph{SIAM Journal on Control and Optimization}, \textbf{47(6)} (2009), 3237--3258.
    
\bibitem{peet_ACC_2010} (MR1124979) [10.2307/2152750]
    \newblock M. Peet and Y. Peet,
    \newblock A parallel-computing solution for optimization of polynomials,
    \newblock \emph{American Control Conference}, (2010).
      
\bibitem{hilbert1} (MR1124979) [10.2307/2152750]
    \newblock D. Hilbert,
    \newblock Uber die Darstellung definiter Formen als Summe von Formen quadratens,
    \newblock \emph{Math. Ann.}, \textbf{32} (1888), 342--350.    
      
\bibitem{hilbert2} (MR1124979) [10.2307/2152750]
    \newblock D. Hilbert,
    \newblock Uber ternare definite Formen,
    \newblock \emph{Acta Math.}, \textbf{17} (1893), 169--197.     
  
\bibitem{motzkin} (MR1124979) [10.2307/2152750]
    \newblock T.S. Motzkin,
    \newblock The arithmetic-geometric inequality,
    \newblock \emph{Symposium on Inequalities, Academic Press}, \textbf{253} (1967), 205--224.

\bibitem{Hil17_reznick} (MR1124979) [10.2307/2152750]
    \newblock B. Reznick,
    \newblock Some concrete aspects of Hilbert's 17th problem,
    \newblock \emph{Contemporary Mathematics}, \textbf{253} (2000), 251--272.

\bibitem{artin} (MR1124979) [10.2307/2152750]
    \newblock E. Artin,
    \newblock Uber die Zerlegung definiter Funktionen in Quadra,
    \newblock \emph{Quadrate, Abh. Math. Sem. Univ. Hamburg}, \textbf{5} (1927), 85--99.

\bibitem{polya_book} (MR1124979) [10.2307/2152750]
    \newblock G. Hardy and J. E. Littlewood and G. P{\'o}lya,
    \newblock Inequalities,
    \newblock \emph{Cambridge University Press}, \textbf{5} (1934), 85--99.
    
\bibitem{habicht} (MR1124979) [10.2307/2152750]
    \newblock W. Habicht,
    \newblock Uber die Zerlegung strikte definiter Formen in Quadrate,
    \newblock \emph{Commentarii Mathematici Helvetici}, \textbf{12(1)} (1939), 317--322.    
  
\bibitem{reznick_no_denominator} (MR1124979) [10.2307/2152750]
    \newblock B. Reznick,
    \newblock On the absence of uniform denominators in Hilbert’s 17th problem,
    \newblock \emph{Proceedings of the American Mathematical Society}, \textbf{133(10)} (2005), 2829--2834.
    
\bibitem{bernstein_1915} (MR1124979) [10.2307/2152750]
    \newblock S. Bernstein,
    \newblock Sur la repr sentation des polynomes positif,
    \newblock \emph{Soobshch. Har'k. Mat. Obshch.}, \textbf{2(14)} (1915), 227--228.

\bibitem{positive_on_interval_reznick_powers} (MR1124979) [10.2307/2152750]
    \newblock V. Powers and B. Reznick,
    \newblock Polynomials that are positive on an interval,
    \newblock \emph{Transactions of the American Mathematical Society}, \textbf{352(10)} (2000), 4677--4692.    

\bibitem{handelman} (MR1124979) [10.2307/2152750]
    \newblock D. Handelman,
    \newblock Representing polynomials by positive linear
functions on compact convex polyhedra,
    \newblock \emph{Pac. J. Math}, \textbf{132(1)} (1988), 35--62.   

\bibitem{schweighofer_nonnegativity} (MR1124979) [10.2307/2152750]
    \newblock M. Schweighofer,
    \newblock Certificates for nonnegativity of polynomials with zeros on compact semialgebraic sets,
    \newblock \emph{manuscripta mathematica}, \textbf{117(4)} (2005), 407--428.
  
\bibitem{scheiderer_survey} (MR1124979) [10.2307/2152750]
     \newblock C. Scheiderer,
     \newblock \emph{Positivity and sums of squares: a guide to recent results},
     \newblock Emerging applications of algebraic geometry, Springer New York, 2009, 271--324.    
    
\bibitem{laurent_survey} (MR1124979) [10.2307/2152750]
     \newblock M. Laurent,
     \newblock \emph{Sums of squares, moment matrices and optimization over polynomials},
     \newblock Emerging applications of algebraic geometry, Springer New York, 2009, 157--270. 
      
\bibitem{delzell_book} (MR1124979) [10.2307/2152750]
     \newblock A. Prestel and C. Delzell,
     \newblock \emph{Positive polynomials: from Hilbert’s 17th problem
to real algebra},
     \newblock Springer New York, 2004. 
     
\bibitem{reznick_powers_polyhedra} (MR1124979) [10.2307/2152750]
    \newblock V. Powers and B. Reznick,
    \newblock A new bound for Polya's Theorem with applications to polynomials positive on polyhedra,
    \newblock \emph{Journal of Pure and Applied Algebra}, \textbf{164} (2001), 221--229.     
     
\bibitem{SDP_SIAMbook_parrilo} (MR1124979) [10.2307/2152750]
     \newblock G. Blekherman, P. Parrilo and R. Thomas,
     \newblock \emph{Semidefinite optimization and convex algebraic geometry},
     \newblock SIAM Philadelphia, 2013. 

\bibitem{lasserre2001} (MR1124979) [10.2307/2152750]
    \newblock J. Lasserre,
    \newblock Global Optimization with Polynomials and the Problem of Moments,
    \newblock \emph{SIAM Journal on Optimization}, \textbf{11(3)} (2001), 796--817.  
    
\bibitem{khalil} (MR1124979) [10.2307/2152750]
     \newblock H. Khalil,
     \newblock \emph{Nonlinear Systems},
     \newblock Prentice Hall; Third edition, New Jersey, 2002. 
     
\bibitem{bliman_poly_solution_2004} (MR1124979) [10.2307/2152750]
    \newblock P. Bliman,
    \newblock An existence result for polynomial solutions of parameter-dependent LMIs,
    \newblock \emph{Systems and Control Letters}, \textbf{51(3)} (2004), 165--169.  

\bibitem{gahinet_chilali_TAC1996} (MR1124979) [10.2307/2152750]
    \newblock P. Gahinet, P. Apkarian and M. Chilali,
    \newblock Affine parameter-dependent Lyapunov functions and real parametric uncertainty,
    \newblock \emph{IEEE Transactions on Automatic Control}, \textbf{41(3)} (1996), 436--442

\bibitem{oliveira_peres_2005} (MR1124979) [10.2307/2152750]
    \newblock R. Oliveira, P. Peres,
    \newblock Stability of polytopes of matrices via affine parameter-dependent Lyapunov functions: Asymptotically exact LMI conditions,
    \newblock \emph{Linear Algebra and its Applications}, \textbf{405} (2005), 209--228.

\bibitem{bliman_homogeneous} (MR1124979) [10.2307/2152750]
    \newblock P. Bliman, R. Oliveira, V. Montagner and P. Peres,
    \newblock Existence of Homogeneous Polynomial Solutions for
Parameter-Dependent Linear Matrix Inequalities with Parameters in
the Simplex,
    \newblock \emph{IEEE Conference on Decision and Controls}, (2006).

\bibitem{reza_peet_acc_2012} (MR1124979) [10.2307/2152750]
    \newblock R. Kamyar and M. Peet,
    \newblock Existence of Decentralized computation for robust stability analysis of large state-space systems using Polya's theorem,
    \newblock \emph{American Control Conference}, (2012).


\bibitem{bliman_SIAM_2004} (MR1124979) [10.2307/2152750]
    \newblock P Bliman,
    \newblock A convex approach to robust stability for linear systems with uncertain scalar parameters,
    \newblock \emph{SIAM journal on Control and Optimization}, \textbf{42(6)} (2004), 2016--2042.
    
\bibitem{chesi_hypercube_2005} (MR1124979) [10.2307/2152750]
    \newblock G. Chesi,
    \newblock Establishing stability and instability of matrix hypercubes,
    \newblock \emph{System and Control Letters}, \textbf{54} (2005), 381--388.

\bibitem{peet_papa_TAC2012} (MR1124979) [10.2307/2152750]
    \newblock M. Peet and A. Papachristodoulou,
    \newblock A converse sum of squares Lyapunov result with a degree bound,
    \newblock \emph{IEEE Transactions on Automatic Control}, \textbf{57(9)} (2012), 2281--2293.
    
\bibitem{sigurdur_existence_2012} (MR1124979) [10.2307/2152750]
    \newblock P. Giesl and S. Hafstein,
    \newblock Existence of piecewise linear Lyapunov functions in arbitary
dimensions,
    \newblock \emph{Discrete and Continuous Dynamical Systems}, \textbf{32} (2012), 3539--3565.   
    
\bibitem{chesi_2005} (MR1124979) [10.2307/2152750]
    \newblock G. Chesi, A. Garulli, A. Tesi and A. Vicino,
    \newblock LMI-based computation of optimal quadratic Lyapunov functions for odd polynomial systems,
    \newblock \emph{International Journal of Robust and Nonlinear Control}, \textbf{15(1)} (2005), 35--49. 
        
\bibitem{Gahinet_1994_LMI_hinf} (MR1124979) [10.2307/2152750]
    \newblock P. Gahinet, P. Apkarian,
    \newblock A linear matrix inequality approach to H infinity control,
    \newblock \emph{International Journal of Robust and Nonlinear Control}, \textbf{4(4)} (1994), 421--448.    
    


\end{thebibliography}
\end{document}